\input amstex
\documentstyle{amsppt}
\magnification=\magstep1
\nologo
\frenchspacing
\def\ba#1{{\bar #1}}
\def\bo#1{{\bold #1}}
\def\ca#1{{\Cal #1}}
\def\go#1{{\goth #1}}

\def\back{{\backslash}}
\def\bH{\overline{\bold H}}
\def\cosh{{\text{\rm cosh}}}
\def\congr{ \equiv }
\def\disc{\text{disc\,}}
\def\Gal{{\text{\rm Gal}}}
\def\gcd{{\text{\rm gcd}}}
\def\GL{{\text{\rm GL}}}
\def\half{{\textstyle{1\over 2}}}
\def\H{\bold H}

\def\iso{ \cong }
\def\id{{\text{\rm id}}}
\def\Im{{\text{\rm Im}}}

\def\legendre{\overwithdelims()}

\def\mod{\bmod\,}

\def\ma#1#2#3#4{\left({{#1}\atop{#3}}\; {{#2}\atop{#4}}\right)}
\def\pma#1#2#3#4{\pmatrix {#1}&{#2}\cr {#3}&{#4}\cr\endpmatrix}
\def\N{{\text{\rm N}}}
\def\PSL{{\text{\rm PSL}}}
\def\Re{{\text{\rm Re}}}

\def\sumprime_#1 {\setbox0=\hbox{$\scriptstyle{#1}$}
   \setbox2 =\hbox{$\displaystyle{\sum}$}
   \setbox4 =\hbox{${}'\mathsurround=0pt$}
   \dimen0=.5\wd0  \advance\dimen0 by-.5\wd2
   \ifdim\dimen0>0pt
   \ifdim\dimen0>\wd4 \kern\wd4 \else\kern\dimen0\fi\fi
   \mathop{{\sum}'}_{\kern-\wd4 #1}}
\def\SL{{\text{\rm SL}}}
\def\sign{{\text{\rm sign}}}
\def\Tr{{\text{\rm Tr}}}
\newcount\refCount
\def\newref#1 {\advance\refCount by 1
\expandafter\edef\csname#1\endcsname{\the\refCount}}
\newref CO     
\newref DO   
\newref GE    
\newref GR    
\newref GZe 
\newref GZd  
\newref GZt  
\newref HEc  
\newref HEj   
\newref IW     
\newref LA     
\newref MI     
\newref NE     
\newref RA    
\newref SI     
\newref YZ    
\newref ZA     

\def\today{\ifcase\month\or January \or February\or
March\or April\or May\or June\or July\or August\or
September\or October\or November\or December\fi
\space\number\day, \number\year}

\newcount\thmnum
\thmnum=1
\def\thmnumber{\the\thmnum. \global\advance\thmnum  by 1}

\newcount\eqnum
\eqnum=1
\def\eqnumber{\the\eqnum \global\advance\eqnum  by 1}

\def\eqn#1{{\advance\eqnum by -#1\negthinspace (\the\eqnum)}}
\def\thmn#1{{\advance\thmnum by -#1 \the\thmnum}}

\topmatter
\title
On singular moduli for level 2 and 3
\endtitle
\author
Hans Roskam
\endauthor
\abstract
Gross and Zagier proved a formula for the 
absolute norm $N(j(\alpha_1)-j(\alpha_2))$ of
a difference of singular values of the modular
function $j$.
We formulate and prove the analogues of their result
for a number of functions of level 2 and 3.
\endabstract
\endtopmatter
\voffset-.8in 
\document
\head 1. Introduction
\endhead
\noindent
In a famous paper [\GZe], Gross and Zagier established 
an explicit formula for the expression
$$
J(d_1, d_2)= \Bigl(\prod_{[\tau_1], [\tau_2]\atop \disc\tau_i=d_i}
             \big (j(\tau_1)-j(\tau_2)\big)\Bigr)^{4\over w_1w_2}.
$$
Here $j$ is the elliptic modular function on the complex upper half plane
$\H=\{z\in\bo C:\Im(z)>0\}$,
$d_1$ and $d_2$ are negative coprime fundamental quadratic discriminants,
$w_1$ and $w_2$ are the number of roots of unity in the corresponding
imaginary quadratic orders,
and $[\tau_i]$ denotes the equivalence class of $\tau_i\in\H$ under the
natural action of $\text{SL}_2(\bo Z)$.
By $\disc\tau_i=d_i$ we mean that $\tau_i$ is imaginary quadratic,
and that its irreducible polynomial over~$\bo Z$ has discriminant~$d_i$.
We know by the theory of complex multiplication that $J(d_1,d_2)$ is  the ${4\over w_1w_2}$-power of the norm of the algebraic integer $j(\tau_1) - j(\tau_2)$.
In particular if $d_1$ and $d_2$ are both less than $-4$ then $J(d_1,d_2) $ is an integer.

In order to give the Gross-Zagier formula, we define
for primes $p$ satisfying ${d_1d_2\legendre p}\not=-1$
$$\varepsilon(p)=\cases
 {d_1\legendre p}&\hbox{if\ } \gcd(p,d_1)=1;\cr
 {d_2\legendre p}&\hbox{if\ } \gcd(p,d_2)=1.\cr
\endcases$$
We extend $\varepsilon$ multiplicatively to products
$m=\prod_i p_i^{a_i}$ of such prime numbers, and put
$$
F(m)=\prod_{\scriptstyle nn ^\prime =m\atop\scriptstyle n,n ^\prime>0}n^{\varepsilon(n ^\prime)}
$$
for such a product.
Furthermore we set $F(m) =1 $ for all $m\in\bo Q$ that are not of the above form.
\bigskip\noindent\bf Theorem\,\thmnumber\sl
Using the above notation and setting $D =d_1d_2 $  the following formula holds:
$$
J(d_1,d_2)^2  =\pm \prod_{\scriptstyle x\in\bo Z\atop\scriptstyle  x^2<D}
\hskip-0.1cm F\Big( {D-x^2\over 4} \Big).
$$
\medskip \noindent\rm 
The above product can be restricted to those $x\in\bo Z $ that satisfy \hbox{$x ^2\congr D\mod 4$.}
For these integers one proves the equality $\varepsilon({D- x ^2\over4}) = -1$ using quadratic reciprocity.
Furthermore one can prove that  for any positive integer $m$ with $\varepsilon(m)=-1$ the value
 $F(m)$ is a prime power [\CO, page 306--307].
More precisely we have
$$
F(m)=\cases
p^{(a+1)\prod_{j=1}^s(b_j+1)}&\vtop{\hsize=7cm \noindent
if $m=p^{2a+1}\cdot\prod_{i=1}^rp_i^{2a_i}\cdot\prod_{j=1}^sq_j^{b_j}$,
where
$p,p_i$ and $q_j$ are distinct primes satisfying
$\varepsilon(p)=\varepsilon(p_i)=-1$ and $\varepsilon(q_j)=1$;}\cr
\noalign{\medskip}
1&\hbox{otherwise.}\cr
\endcases
$$
Using this formula the following corollary is immediate.
\medskip\noindent\bf Corollary\,\thmnumber\sl
If $p$ is a prime dividing $J(d_1,d_2) ^2 $ then ${d_1\legendre p}\ne 1, {d_2\legendre p}\ne 1$ and $p$ divides a positive integer of the form
${D-x^2\over 4}$. In particular, $p\le{D\over 4}$.
\bigskip \noindent\rm For algebraic numbers $\beta_1,\beta_2 $ 
we set 
$$
\N(\beta_1,\beta_2)=|\N_{\bo Q(\beta_1,\beta_2)/\bo Q}(\beta_1-\beta_2)|.
$$
Hence the above theorem gives an expression for  $\N\big (j(\alpha_1), j(\alpha_2)\big)^{8\over w_1w_2}$,
where $\alpha_i\in\bo H $ are arbitrarily under the condition $\disc \alpha_i=d_i$. 
We will prove similar expressions for $\N\big (f(\alpha_1),f(\alpha_2)\big)$
for a number of  modular functions $f$ of level~2 and~3.
Here $\alpha_i\in\bo H $ of discriminant $d_i$ is chosen such that $f(\alpha_i)$ is a class invariant, 
i.e.  such that $j(\alpha_i)$ and $f(\alpha_i)$  generate the same field over $\bo Q(\sqrt{d_i})$.
 
Recall the definition of the $j$-function 
$$
j(z)= 12 ^3{g_2(z)^3\over \Delta(z) } = 12 ^3+6 ^ 6 {g_3(z)^2\over \Delta(z) } 
$$
where $g_2(z),g_3(z)$ and $\Delta(z)$  are the modular forms on $\SL_2(\bo Z)$  of weight 4, 6 and~12 respectively defined by
$$\eqalign{
g_2(z)=60\hskip-0.4cm\sum_{\scriptstyle m,n\in\bo Z\atop\scriptstyle 
(m,n)\ne(0,0)}{1\over(mz+n)^4},\
\
&
\ \ g_3(z)=140\hskip-0.4cm\sum_{\scriptstyle m,n\in\bo Z\atop\scriptstyle 
(m,n)\ne(0,0)}{1\over(mz+n)^6}
\cr}$$
and $\Delta (z) = g_2 ^3(z)-27g_3 ^2(z)=(2\pi)^{12}\eta(z)^{24}$. 
Here $\eta(z)$ is the Dedekind eta function
$$
\eta(z)=q^{1\over 24}\prod_{n=1}^\infty(1-q^n)
$$
with $q =q(z) =e^{2\pi iz}$.

The first two functions that we will study are the Weber functions
$$\eqalign{
\gamma_3(z)=6^3{g_3(z)\over{(2\pi)}^6\eta(z)^{12}},\ \ &\ \ \gamma_2(z)=12\,{g_2(z)\over{(2\pi)}^4\eta(z)^{8}},\cr
}$$
which are modular functions of level 2 and 3 respectively.
The relations
$$\eqalign{
\gamma_3(z)^2=j(z)-12 ^3,\ \ &\ \ \gamma_2(z)^3=j(z)\cr
}\leqno(\eqnumber)$$
imply that for $f\in\{\gamma_3,\gamma_2\} $ any prime divisor of $\N\big(f(\alpha_1),f(\alpha_2)\big) $ also divides
 $\N\big(j(\alpha_1),j(\alpha_2)\big)$.
In fact for $f =\gamma_3$ the converse also holds.
Before we prove this we fix the following notation which we use throughout this paper:
\bigskip 
\centerline{\vbox{\hsize=11cm\sl\noindent  The discriminants $d_1$ and $d_2$ of the imaginary quadratic fields $K_1$ and $K_2$ are relatively prime and $w_1,w_2 $ and $h_1, h_2$ denote the number of roots of unity and the class numbers of their ring of integers. Furthermore we set $D =d_1d_2$ and $h_i ^\prime ={2\over w_i}h_i$ for $i =1,2 $.}}
\bigskip\noindent\bf Theorem\,\thmnumber\sl
For $\alpha_1,\alpha_2\in\H$ of discriminant $d_1$ and $d_2$ respectively,
the following formula holds:
$$
\N\big(\gamma_3(\alpha_1),\gamma_3(\alpha_2)\big)=\N\big (j(\alpha_1),j(\alpha_2)\big)^e
$$
with $e =2 $ if neither $d_1$ nor $d_2$ is equal to $-4$ and $e =1 $ otherwise. 
\bigskip \noindent\bf Proof.
\rm First assume that neither $d_1$ nor $d_2$ is equal to $-4$.
We claim that in this case $\gamma_3(\alpha_i) $ generates a quadratic extension of $\bo Q(j(\alpha_i))$.
Namely, if $d_i$ is odd then $\gamma_3(\alpha_i)$ is conjugate over $\bo Q$ to one of the numbers
$\pm \gamma_3({-1+\sqrt{d_i}\over2})$.
As $j({-1+\sqrt{d_i}\over2})$ is real and less then $12 ^3 $  the claim follows from the first equality in \eqn1.
If $d_i$ is even and different from $-4$,
then $K_i(\gamma_3(\alpha_i)) $ is the ray class field of conductor 2 of $K_i$ which is quadratic over
 $K_i(j(\alpha_i))$  [\GE, theorem 10].
We conclude that if $\disc\alpha_i\ne -4$ then $\gamma_3(\alpha_i)$ is of degree $2h_i $ over $\bo Q$
and hence has $\{\gamma_3(\tau_i), -\gamma_3(\tau_i)\}_{[\tau_i],\disc\tau_i =d_i}$ as a complete set of conjugates over $\bo Q$.

As $K_i (\gamma_3(\alpha_i))$ is the Hilbert class field or the ray class field of conductor 2 of $K_i$, depending on whether $d_i$ is odd or even, its subfield $\bo Q (\gamma_3(\alpha_i)) $  is unramified at the primes 
not dividing $d_i$.
By assumption $d_1$ and $d_2$ are relatively prime.
The fields $\bo Q(\gamma_3(\alpha_1)) $ and $\bo Q(\gamma_3(\alpha_2))$ are therefore linearly disjoint over $\bo Q$.
Systematically using
$$(x^2-y^2)^2=(x-y)(x+y)(-x+y)(-x-y)$$
we find
$\N\big(\gamma_3(\alpha_1),\gamma_3(\alpha_2)\big)=\N\big (j(\alpha_1),j(\alpha_2)\big)^2$.

If $\alpha_i$ is of discriminant $-4$ then $\gamma_3(\alpha_i) =0 $ and the formula for the absolute norm of $\gamma_3(\alpha_1) -\gamma_3(\alpha_2)$ follows easily. 
\hfill$\qed$ 
\bigskip\noindent Our result for the function $f =\gamma_2$ is less complete.
As mentioned above we normalise $\alpha_i$ such that $\gamma_2(\alpha_i)$ is a class invariant.
Although this normalisation works for all discriminants $d_i$ that are coprime to 3,
we only obtained a formula in the case $d_1\congr d_2\congr 2\mod 3$. 
 \bigskip\noindent\bf Theorem\,\thmnumber\sl
Assume that $d_1$ and $d_2$ are both congruent to $2 $ modulo $3 $.
For $i=1,2$ define $\alpha_i\in\H$ by
$$\alpha_i=
\cases
{{-3+\sqrt{d_i}}\over 2}&\hbox{if\ }d_i\congr 1\mod 4;\cr
{\sqrt{d_i}\over 2}&\hbox{if\ }d_i\congr 0\mod 4.\cr
\endcases
$$
Then the following formula holds:
$$
\N\big(\gamma_2(\alpha_1),\gamma_2(\alpha_2)\big)^{8\over{w_1w_2}}=3^{6h_1^\prime h_2^\prime}\cdot\hskip -0.1cm
\prod_{\scriptstyle x\in\bo Z\atop\scriptstyle x ^2 <D}\hskip -0.1cm F\Big({{D-x^2}\over 36}\Big).
$$
\rm\smallskip \noindent As a corollary we find that the primes $p\ne 3$  dividing $\N\big(\gamma_2(\alpha_1),\gamma_2(\alpha_2)\big)$ satisfy $p<{D\over 36}$.
If both $d_1$ and $d_2$ are not equal to $-4$, the exponent ${8\over w_1w_2} $
is equal to~2.  
In this case one obtains a formula for $\N\big(\gamma_2(\alpha_1),\gamma_2(\alpha_2)\big) $ 
by restricting the above product to positive $x$ and by replacing the exponent $6h_1^\prime h_2^\prime$ by $3h_1h_2$.

Next we consider the 24-th powers of the classical Weber $\go f$-functions:
$$\eqalign{
\omega(z)&={\Delta({z+1\over 2})\over\Delta(z)}= -q^{-\half}\prod_{n=1}^\infty(1+q^{n-\half})^{24},\cr
\omega_1(z)&={\Delta({z\over 2})\over\Delta(z)}=q^{-\half} \prod_{n=1}^\infty(1-q^{n-\half})^{24},\cr
\omega_2(z)&=2^{12}{\Delta({2z})\over\Delta(z)}=2^{12}q\prod_{n=1}^\infty(1+q^{n})^{24},\cr
}$$
where for  $z\in\bo H $ and $a\in\bo Q$ we set $q^a=e^{2\pi i a z}$.
These functions are modular of level 2 and satisfy the polynomial equation [\CO, theorem 12.17]
$$(X-\omega)(X-\omega_1)(X-\omega_2)=(X+16)^3-jX.\leqno(\eqnumber)$$
\bigskip\noindent\bf Theorem\,\thmnumber\sl
Assume that $d_1$ and $d_2$ are both congruent to $1 $ modulo $8$.
 For $i=1,2$ define $\alpha_i\in\H$ by $\alpha_i={{-1+\sqrt{d_i}}\over2}$.
 Then the following formulas hold:
 $$\displaylines{
 \N\big(\omega(\alpha_1),\omega(\alpha_2)\big)=\N\big(\omega_1(\alpha_1),\omega_1(\alpha_2)\big)=
 2^{12h_1h_2}\cdot\hskip-0.1cm
 \prod_{\scriptstyle x\in\bo Z\atop\scriptstyle  x^2<D}\hskip -0.1cm
 F\Big({{D-x^2}\over 8}\Big),\cr
 \noalign{\medskip}
 \N\big(\omega_2(\alpha_1),\omega_2(\alpha_2)\big)=\hskip -0.2 cm\prod_{\scriptstyle x\in\bo Z\atop\scriptstyle  0<x<\sqrt{D}}\hskip-0.3cm F\Big({{D-x^2}\over 16}\Big).\cr
 }$$
\rm\medskip \noindent 
The above formula for $\omega_2$ was conjectured by Yui and Zagier in [\YZ, formula $5_?$].

In the table below we have listed the prime factorization of $\N\big (f(\alpha_1) -f(\alpha_2)\big)$ for each of the functions $f\in\{j,\gamma_2,\omega,\omega_2\} $ and for two choices of the pair $(\alpha_1,\alpha_2)$. 
Here $\alpha_i\in\bo H $ is of discriminant $d_i$ and normalised depending on $f$ as in the theorems above.
The class numbers of the quadratic fields in the first column are equal to $h_1 = 1 $ and $h_2 = 4 $, in the second column they equal $h_1 =3 $ and $ h_2 =7$. 
\def\mark{\hbox to 4pt{\hss$\cdot$\hss}}
\vskip0.35cm \noindent
{
\hfill \vbox{\halign{ $#$\quad &$#$\quad &$#$ \cr
&\hfill (d_1,d_2) = (-7, -55)\hfill  &\hfill (d_1,d_2) = (-31, -151)\hfill  \cr
j & \hfill 3 ^ {26}  5 ^6  19 ^ 3  47 ^2 &\hfill 3^ {140}   13^ {21}   23^ {11}  53^6   61^2   73^5   79^4   83^4 89^2   179^2   449^2   557^2 \cr
\gamma_2 &\hfill  3 ^ {14}  5 ^ 2 &\hfill 3^ {77}  13^7   23^5   61^2\cr
\omega &\hfill  2 ^ {48}  3 ^ {34}  5 ^ 8  19 ^4  47^2 &\hfill  2^ {252}   3^ {190}   13^ {28}   23^{12}   53^8   61^2   73^8   79^4   83 ^6   89^2 179^2   449^2  557^2\cr
\omega_2 &\hfill  3 ^ 8  5 ^ 2  19 &\hfill  3^ {50}   13^7   23\mark  53^2   73^3   83^2\cr
}}\hfill }
\vskip0.35cm  \noindent Gross and Zagier gave two proofs of theorem \thmn5, one algebraic and one analytic. 
Their algebraic proof uses the reduction theory of elliptic curves and was restricted to the case of prime discriminants. Later Dorman extended it to the general case [\DO].
Our proof of theorems \thmn2 and \thmn1 is an adaptation of the analytic proof of Gross and Zagier.
It consists of three steps, sections 2, 3 and 4, which are more or less independent.
In the outline below we concentrate on the function~$\gamma_2$.

In section 2 we use Shimura's reciprocity law to compute the conjugates of $\gamma_2(\alpha)$ where $\alpha\in\bo H$ is of discriminant $d\congr 2 \mod 3$ and normalised as in theorem~\thmn2.
These conjugates can be written in the form $\gamma_2(\tau)$ 
with $\tau $ in some finite set  $S_d$ of elements $\tau\in\bo H$  with $\disc\tau =d $ and $\Tr_{\bo Q (\sqrt{d})/\bo Q}(\tau)\congr 0\mod 3$.
Next in section~3 we characterise $h(z_1,z_2) =\log |\gamma_2(z_1) -\gamma_2(z_2) | $ as the unique
symmetric and harmonic function on $\bo H\times\bo H $ with certain invariance and growth conditions. 
Using Legendre functions and Eisenstein series we then build a function satisfying the same properties.
By summing this function over  $(z_1,z_2)\in S_{d_1}\times S_{d_2} $  we arrive at a complicated looking expression for $\log \N\big (\gamma_2(\alpha_1),\gamma_2(\alpha_2)\big)$ (theorem  \thmn{-7} below).
These kinds of expressions were recognised by Gross and Zagier as being related to the Fourier coefficients  of certain holomorphic modular forms of weight~2 on congruence subgroups of $\SL_2(\bo Z)$.
To make this precise we study in section~4 a family of non-holomorphic Hilbert modular forms and,
via restriction to the diagonal and holomorphic projection, the corresponding family of holomorphic modular forms.
Finally in section 5 we conclude that $\log \N\big (\gamma_2(\alpha_1),\gamma_2(\alpha_2)\big)$ is up to
a simple expression equal to the first Fourier coefficient of one of these holomorphic modular forms.
As the corresponding congruence subgroup has genus zero, this Fourier coefficient vanishes and the proof of theorem \thmn2 is complete.

\head 2. Computing conjugates
\endhead
\noindent 
An element $M\in\PSL_2(\bo Z) $ represented by the matrix $\ma abcd\in\SL_2(\bo Z) $  acts on the upper half plane $\bo H$  by the linear fractional transformation
$$\textstyle z \mapsto Mz={az +b\over cz +d}.  $$
In the sequel we will identify the transformation $M$ with the matrix $\ma abcd$ and write $M =\ma abcd $
instead of $M \congr \ma abcd \mod \{\pm\id\}$.
The left  $\PSL_2(\bo Z)$-action on $\bo H$  induces a right action on functions $f$ on $\bo H$ by $(f\circ M)(z)=f(Mz)$.
Fix a positive integer $N$.
The (projective) principal congruence modular group $\overline{\Gamma}(N)$ is defined as the kernel of the map
 $\PSL_2(\bo Z)\rightarrow\SL_2 (\bo Z/N\bo Z)/\{\pm\id\} $ induced by the natural map $\bo Z\rightarrow\bo Z/N\bo Z$. 
A modular function of level $N$ is a meromorphic function  $f$ on~$\bo H$ that is invariant under $\overline{\Gamma}(N)$,
 i.e.  $f\circ M =f$ for all $M\in\overline{\Gamma}(N)$, and that is `meromorphic at the cusps'.  
Recall that the cusps are the points $x\in\bo P ^1(\bo Q)$ which are on the boundary of $\bo H$.
To clarify the condition `meromorphic at the cusps' we need the $N$-th root $q^{1\over N} =e^{2\pi iz\over N}$ of the function $q =q(z) =e^{2\pi iz}$. 
A $\overline{\Gamma}(N)$-invariant function $f$ satisfies $f(z+N) =f(z) $ because $\ma 1N01\in\overline{\Gamma}(N)$.
Hence there is a meromorphic function $f ^* $ on the punctured disk $\{q ^{1\over N}:0 < |q^{1\over N} | < 1\}$  
such that $f(z) =f ^*(q^{1\over N})$.
If~$f ^* $ is meromorphic at $q^{1\over N} =0$  we say that $f$ is meromorphic at the cusp $\infty$.
The $q$-expansion of $f$ is by definition the Laurent expansion of~$f ^*$ at $q^{1\over N} =0$.
Finally we say that $f$ is meromorphic at the cusps if $f\circ M$ is meromorphic at $\infty$ for all
 $M\in\PSL_2(\bo Z)$.
Note that the behaviour of $f$ near the cusp $x\in\bo P ^1(\bo Q)$ is reflected by the behaviour of $f\circ M$
near $\infty$ if $M x =\infty$.   

An equivalent way of introducing modular functions of level $N$ is the following.
Let $\bH =\bo H\cup\bo P ^1 (\bo Q)$ be the extended upper half plane and extend the $\PSL_2(\bo Z)$-action
on $\bo H$ in the obvious way to $\bH$.
The orbit space $\overline{\Gamma}(N)\back\bH$ can be given a complex structure [\MI, \S 1.8].
The addition of the set of cusps $\bo P ^1(\bo Q)$ to $\bo H$ causes the resulting Riemann surface to be compact.
The meromorphic functions on this Riemann surface correspond to the modular functions of level $N$.

Let $F_N$ be the field of modular functions of level $N$ for which the $q$-expansion is rational over 
$\bo Q(\zeta_N)$, with $\zeta_N =e^{2\pi i\over N}$ a $N$-th root of unity.
Fix an imaginary quadratic field $K\subset\bo C$ with ring of integers $\ca O$. 
The evaluation of $f\in F_N $ at $\theta\in\bo H\cap K$ is called a singular modulus.
By the theory of complex multiplication a singular modulus generates an abelian extension of the quadratic number field~$K$.
More precisely the field generated over $K$ by $f(\theta)$ with fixed $\theta$ as above and $f$ ranging over those functions in $F_N$ that are defined at $\theta$, is equal to the ray class field of conductor  $N$ [\LA, page 128]. 

For example in the case $N =1 $  we have $F_1=\bo Q(j) $. 
If  $\theta\in\bo H$ generates the ring of integers $\ca O$ the field $K(j(\theta))$ is equal to the Hilbert class field $H$ of $K$, the maximal abelian unramified extension of $K$. 
By class field theory the Artin map supplies an isomorphism between the ideal class group $C$ of $K$ and the Galois group $\Gal(H/K)$.
To describe the action of this group on $j(\theta)$ we represent elements of $C$ by  $\PSL_2(\bo Z)$-equivalence  classes of primitive positive definite quadratic forms of discriminant equal to the discriminant of $K$.
For $d\congr 0,1\mod 4$ a negative integer let
$$
\ca Q_d=\{[a,b,c] : a,b,c\in\bo Z, a>0, b^2-4ac=d\}
$$
be the set of positive definite quadratic forms of discriminant $d$. 
If $d$ is the discriminant of $K$ and $[a, -b,c]\in\ca Q_d $ we have
$$\textstyle
j(\theta)^{[a, -b,c]} =j({-b+\sqrt{d}\over2a}).\leqno(\eqnumber)
$$
Because the function $j$ is $\PSL_2(\bo Z)$-invariant, the above value only depends on the $\PSL_2(\bo Z)$-equivalence class of the quadratic form $[a, -b,c]$.
Formula \eqn1  is a reformulation of the well known formula $(\go b,K/\bo Q)j(\ca O) =j({\go b}^{-1}) $,  where $(\go b,K/\bo Q)$ denotes the Artin automorphism of the  $K$-ideal $\go b$. 

More generally for $f\in F_N $ and $\theta$ a generator of $\ca O$ the singular modulus $f(\theta)$ 
lies in the ray class field of conductor $N$ of $K$, a field containing the Hilbert class field $H$ of $K$.
In the examples treated in this paper $f(\theta)$ does in fact generate $H$ over $K$.
In this case $f(\theta)$ is called a class invariant and we will restrict ourselves to this case in the discussion below.

To obtain the generalization of \eqn1 to class invariants we use Shimura's reciprocity law.
First we need to describe the Galois theory of the modular function fields $F_N$ [\LA, chapter 6].
The field $F_N$ is Galois over $F_1$ with group isomorphic to $\GL_2(\bo Z/N\bo Z)/\{\pm 1\}$.
For $M\in\GL_2(\bo Z/N\bo Z)/\{\pm 1\}$ we write $M=\ma 100\delta \ba M$ with $\delta=\det M\in(\bo Z/N\bo Z) ^* $ and  $\ba M\in\SL_2(\bo Z/N\bo Z)/\{\pm1\}$. 
Lift $\ba M$ to an element ${\hat M}\in\PSL_2(\bo Z)$ and let $\sigma_\delta$ be the automorphism of $\bo Q(\zeta_N)$
induced by $\zeta_N\mapsto\zeta_N ^\delta$.
For $f\in F_N $ with $q$-expansion $\sum_k a_kq ^ {k\over N}\in\bo Q(\zeta_N)((q^{1\over N}))$ the Galois action of $M$ is given by
$$
f ^M ={\tilde f}\circ {\hat M},
$$
where ${\tilde f}\in F_N$ has $q$-expansion $\sum_k\sigma_\delta (a_k)q^{k\over N}$.

Fix  the following generator of  $\ca O$: 
$$
\theta=\cases
{\sqrt{d}\over 2} &\text{\rm  if\ } 2\mid d;\cr
{-1+\sqrt{d}\over 2} &\text{\rm  if\ } 2\nmid d.\cr
\endcases
$$
To describe the Galois action on a class invariant $f(\theta)$ for some $f\in F_N $ 
we restrict to the case $N =p $ is prime.
Shimura's reciprocity law [\GE, theorem 20] states that for $f\in F_p$ such that $f(\theta)\in K(j(\theta))$ and $[a,-b,c]\in C$ there exists an element
$M=M(a,b,c)\in\GL_2(\bo Z/p\bo Z)/\{\pm\id\} $ such that
$$
\textstyle{
f(\theta)^{[a,-b,c]}=f^M({-b+\sqrt{d}\over 2a}).
}\leqno(\eqnumber)
$$
In case $a$ is prime to $p$ the element $M$ is represented by the following matrix
$$M=\cases
\pma a{b\over 2}01 &\text{\rm  if\ } 2\mid d;\cr
\pma a{b-1\over 2}01 &\text{\rm  if\ } 2\nmid d.\cr
\endcases\leqno(\eqnumber)$$
As each quadratic form is $\PSL_2(\bo Z)$-equivalent to a form $[a, -b,c]$  with $p\nmid a$ [\CO, lemma 2.3 and 2.25], the above restriction on $a$ is not a serious one.
\bigskip \noindent Before we apply Shimura's reciprocity law to the function $\gamma_2$ we first determine
its stabilizer inside $\PSL_2(\bo Z)$, which is classically denoted by $\Gamma^3$. 
Recall that $S =\ma 0{-1}10$ and $T =\ma 1101$ generate the group $\PSL_2(\bo Z)$.
Their action on the Dedekind eta function is given by
$$\eqalign{
\eta(\textstyle{-1\over z})=\sqrt{-iz}\,\eta(z),\quad&\quad 
\eta(z+1)=\zeta_{24}\,\eta(z)\cr
}\leqno(\eqnumber)$$
where the square root is positive on the positive real axis [\LA, page 253].
Using the definition of $\gamma_2 $ in the introduction we find
$$\eqalign{
\gamma_2\circ S=\gamma_2,\quad&\quad\gamma_2\circ T=\zeta_3^{-1}\gamma_2.
}$$
Hence $\Gamma^3$ is the kernel of the character $\PSL_2(\bo Z)\rightarrow\langle\zeta_3\rangle$ sending the transformation $M$ to $\gamma_2\circ M\over \gamma_2$, and therefore normal of index 3 inside $\PSL_2(\bo Z)$. 
This characterization can be used to show the following equalities [\RA, page 15-19]
$$\eqalign{\textstyle
\Gamma^3
&\textstyle=\left\langle \ma 0{-1}10 ,\ma {-1}{-2}11,\ma {-1}{-1}21 \right\rangle\cr
&\textstyle=\left\{\ma abcd\in\PSL_2(\bo Z): ab +cd\congr 0\mod 3\right\}.\cr
}$$
It follows from the second description that $\overline{\Gamma}(3)$ is contained in $\Gamma^3$, hence
$\gamma_2$ is a modular function of level 3.
By applying the Hurwitz formula [\MI, theorem 4.2.11] we find that the Riemann surface $\Gamma^3\back\bH$ is of genus zero.

For a negative integer $d\congr 0, 1 \mod 4$ we define the  $\PSL_2(\bo Z)$-set 
$$
\ca P_{d}=\{\tau\in\H:a\tau^2+b\tau+c=0 \hbox{ for some $[a,b,c]\in\ca Q_d$}\}.\leqno(\eqnumber)
$$
The map sending $[a,b,c]$ to ${-b+\sqrt{d}\over 2a}\in\H$ is a bijection from $\ca Q_d$ to $\ca P_d$ and we denote the quadratic form corresponding to $\tau\in\ca P_d$ by  $[a_\tau,b_\tau,c_\tau]$.
In case $d$ is a negative fundamental discriminant it follows from the discussion above that $\{j (\tau):\tau\in\PSL_2(\bo Z)\back\ca P_d\}$ is a transitive
 $\Gal(\overline{\bo Q}/\bo Q)$-set. 
We have the following analogous result for $\gamma_2$. 
\bigskip\noindent\bf Proposition\,\thmnumber\sl
Let $d\congr 2\mod 3$ be a negative fundamental discriminant
and define $\alpha\in\H$ by
$$\alpha =
\cases
{{-3+\sqrt{d}}\over 2}&\hbox{if\ }d\congr 1\mod 4;\cr
{\sqrt{d}\over 2}&\hbox{if\ }d\congr 0\mod 4.\cr
\endcases
$$
\item{a.} The algebraic integer $\gamma_2(\alpha)$ is of degree $h$ over $\bo Q$.
\item{b.} The action of $\Gamma^3$ on $\H$ stabilizes the set
$$\ca P_d^{\gamma_2}=\{\tau\in\H:a\tau^2+b\tau+c=0 \hbox{ for some $[a,b,c]\in\ca Q_d$
with $3\mid b$}\}.$$
The orbit set $\Gamma^3\back\ca P_d^{\gamma_2}$ has cardinality $h$ and
$\{\gamma_2(\tau): \tau\in\Gamma^3\back\ca P_d^{\gamma_2}\}$
is a complete set of conjugates of $\gamma_2(\alpha)$ over $\bo Q$.
\rm\bigskip \noindent 
In the proof of proposition 6 we need the following lemma.
\bigskip\noindent\bf Lemma\,\thmnumber\sl
Let $Y\subset X$ be an inclusion of sets and $H\subset G$ be
an inclusion of groups.
Assume $G$ acts on $X$ in such a way that $Y$ becomes a $H$-set.
Let $I$ be a complete set of left coset representatives of $H\subset G$
and assume 
$$\displaylines{
\rlap{\rm (\eqnumber)}\hfill X=\coprod_{g\in I}gY. \hfill\llap{(\text{disjoint union})}\cr
}$$
Then the inclusion $Y\subset X$ induces a bijection 
$$
H\backslash Y\ {\buildrel\sim\over\longrightarrow}\ G\backslash X.\leqno(\eqnumber)
$$
\bigskip\noindent\bf Proof.
\rm Let $y_1,y_2\in Y$ be $G$-equivalent, say $y_1=ghy_2$ with 
$g\in I$ and $h\in H$.  Then $ghy_2$ is an element of $gY\cap Y$.
As this intersection is non-empty if and only if $g\in H$, we find that 
$y_1$ and $y_2$ are $H$-equivalent and the map \eqn1  is injective.
The fact that \eqn1 is surjective is immediate from  \eqn2. \hfill$\qed$ 
\bigskip\noindent\bf Proof of proposition \thmn2. 
\rm 
As $\gamma_2$ is the cube root of $j$ that is real valued on the positive imaginary axis
its $q$-expansion lies in $\bo Q((q^{1\over 3}))$.
With $\alpha$ as in the proposition we find that $\gamma_2(\alpha)$ is real.
Because $\gamma_2(\alpha)$ generates the Hilbert class field of $K$ over $K$ [\CO, theorem 12.2]
we conclude that $\gamma_2(\alpha)$ is of degree $h$ over $\bo Q$.
Hence the conjugates of $\gamma_2(\alpha)$ over $\bo Q$ coincide with those over $K$.

Let $[a,-b,c]$ be a quadratic form of discriminant $d$. 
Because $d$ is congruent to 2 modulo 3, the integer $a$ is prime to 3 and we can 
use  \eqn6 and \eqn5  to compute the conjugate $\gamma_2(\alpha )^{[a,-b,c]}$.
For odd discriminants we find:
$$\textstyle{
\gamma_2(\alpha )^{[a,-b,c]}=\gamma_2\circ T^{-1}({ -1+\sqrt{d}\over 2})^{[a,-b,c]}=
\gamma_2^M({-b+\sqrt{d}\over 2a}),}$$
with $M=T^{-1}\ma a{{b-1\over 2}}01 =
\ma 100a \ma a{(b-3)/2}0{a^{-1}}
\in\GL_2(\bo Z/3\bo Z)$.
As the Fourier expansion of $\gamma_2$ has rational coefficients, the matrix $\ma 100a $ 
acts trivially on $\gamma_2$.
To calculate the action of $\ma a{b-3\over 2}0{a^{-1}}$,
we use the decomposition 
$$\textstyle {
 \ma a{b-3\over 2}0{a^{-1}}=ST^{-a}ST^{-a}ST^{-ab-a}\mod 3
} $$
from [\GE, lemma 6] and obtain
$$\textstyle{
\gamma_2(\alpha )^{[a,-b,c]}=\gamma_2({-b(1+2a^2)+\sqrt{d}\over 2a}).
}$$
A similar calculation shows that this formula is also valid for even 
discriminants.
Because $a$ is prime to 3, we have $ 1+2a^2\congr 0\mod 3$, hence
the conjugates of $\gamma_2(\alpha)$ can be written in the form $\gamma_2(\tau)$ for some 
$\tau\in\ca P^{\gamma_2}_d$.

To conclude the proof, we use lemma  \thmn1 to show that the inclusion 
$\ca P^{{\gamma_2}}_d\subset\ca P_d$ induces a bijection 
$$
\Gamma^3\backslash\ca P^{{\gamma_2}}_d\ {\buildrel\sim\over\longrightarrow}\ \PSL_2(\bo Z)\backslash\ca P_d.$$
For a quadratic form $[a,b,c]$ of discriminant $b^2-4ac\congr 2\mod 3$ with 
$3|b$, the congruence $a\congr c\mod 3$ holds.
Using this congruence one checks that $\Gamma^3$ acts on~$\ca P^{{\gamma_2}}_d$.
The set $I=\{T^k:k=0,1,2\}$ is a complete set of left coset representatives of 
$\Gamma^3$ in $\PSL_2(\bo Z)$.
For $\tau\in\ca P_d$ and $\tau^\prime=T\tau$ we have $b_{\tau^\prime}=b_\tau-2a_\tau$.
As $a_\tau$ is prime to 3 we find that $T$ changes the residue class of $b_\tau$ 
modulo 3 and hence $\ca P_d$ is the disjoint union of $\ca P^{{\gamma_2}}_d,\ T\ca P^{{\gamma_2}}_d$ and  
$T^2\ca P^{{\gamma_2}}_d$. According to lemma  \thmn1, the map above is a bijection. 
\hfill$\qed$
\bigskip \noindent 
Using the notation of theorem \thmn4, we saw in the proof above that  $\gamma_2(\alpha_1)$ and $\gamma_2(\alpha_2)$  are class invariants.
As in the proof of theorem \thmn5 we conclude
that $\bo Q(\gamma_2(\alpha_1))$ and  $\bo Q(\gamma_2(\alpha_2))$ are linearly disjoint over $\bo Q$.  
Hence we find the following corollary of proposition \thmn2  which we need in section 3: 
$$
\log \N\big (\gamma_2(\alpha_1),\gamma_2(\alpha_2)\big) =
\sum_{i=1,2,\tau_i\in{\Gamma}^3\backslash\ca P_{d_i}^{\gamma_2}}
\log |\gamma_2(\tau_1) -\gamma_2(\tau_2) |.\leqno(\eqnumber)
$$
\medskip \noindent
The transformations \eqn5 of the Dedekind eta function imply that the functions $\omega,\omega_1 $
and $\omega_2$ are permuted under the action of $\PSL_2(\bo Z)$.  More precisely we have 
$$\eqalign{
(\omega,\omega_1,\omega_2)\circ S=(\omega,\omega_2,\omega_1),\ \ \ &\ \ \  (\omega,\omega_1,\omega_2)\circ T=(\omega_1,\omega,\omega_2).\cr
}\leqno(\eqnumber)$$
We conclude that $\omega,\omega_1$ and $\omega_2$ are invariant under $T ^2 $ and $ST ^2S$.  
These transformations generate the congruence subgroup $\overline{\Gamma}(2)$.
Recall that the modular function field  $F_2$ of level 2 is a Galois extension of $F_1 =\bo Q(j) $ with Galois group $\GL_2(\bo Z/2\bo Z)/\{\pm\id\}\iso S_3$, the permutation group on 3 symbols.
As each of the functions $\omega,\omega_1$ and $\omega_2 $ generates a different cubic extension of $\bo Q(j)$ we 
find that $F_2$ is generated over $\bo Q$ by these functions.
The stabilizer of $\omega_2$ inside $\PSL_2(\bo Z)$ contains $T$ and $\overline{\Gamma}(2)$ and is therefore equal to 
$$\textstyle
\overline{\Gamma}_0(2) =\left\{\ma abcd\in\PSL_2(\bo Z):c\congr 0 \mod 2\right\},
$$
an index 3 subgroup of $\PSL_2(\bo Z)$. 
The set of cusp classes $\overline{\Gamma}_0(2)\back\bo P ^1(\bo Q) $ is represented by 0 and $\infty$. 
It follows from the $q$-expansions of $\omega_2$ and $\omega_2\circ S =\omega_1$ that
 $\omega_2$ has a single simple zero at $\infty\in\overline{\Gamma}_0(2)\back\bH$.
Consequently the compact Riemann surface $\overline{\Gamma}_0(2)\back\bH$ is of genus zero and has $\bo C(\omega_2)$ as its function field.
\medskip\noindent\bf Proposition\,\thmnumber\sl
Let $d\congr 1\mod 8$ be a negative fundamental discriminant.
\item{a.} The algebraic integers $\omega({-1+\sqrt{d}\over 2})$ and
$\omega_1({-1+\sqrt{d}\over 2})$ are conjugate and of degree~$2h$ over~$\bo Q$.
The algebraic integer $\omega_2({-1+\sqrt{d}\over 2})$ is of degree $h$
over $\bo Q$.
\item{b.} The action of $\overline{\Gamma}_0(2)$ on $\H$ stabilizes the sets
$$\displaylines{
\ca P_d^\omega=\ca P_d^{\omega_1} =\{\tau\in\H:a\tau^2+b\tau+c=0 \hbox{ for some $[a,b,c]\in\ca Q_d$ 
with
$2\mid a$}\}\cr
\noalign{and}
\ca P_d^{\omega_2}=\{\tau\in\H:a\tau^2+b\tau+c=0 \hbox{ for some $[a,b,c]\in\ca Q_d$
with $2\nmid a$}\}.\cr
}$$
The orbit sets $\overline{\Gamma}_0(2)\back\ca P_d^\omega$ and $\overline{\Gamma}_0(2)\back\ca P_d^{\omega_2}$ have
cardinality $2h$ and $h$ respectively.
For $f=\omega$ or $\omega_2$, the set
$\{\omega_2(\tau): \tau\in\overline{\Gamma}_0(2)\back\ca P_d^f\}$
is a complete set of conjugates of $f({-1+\sqrt{d}\over 2})$ over $\bo Q$.
\bigskip\noindent\bf Proof.
\rm First we study the action of $\overline{\Gamma}_0(2)$ on $\ca P_d$. 
Let $\tau\in \H$ satisfy the quadratic equation $aX^2+bX+c=0$.
For $\ma xyzw \in\PSL_2(\bo Z)$, the element
$\ma xyzw \tau$ satisfies the equation $a^\prime X^2+b^\prime X+c^\prime=0$ 
with 
$$\leqalignno{
a^\prime=&\ aw^2-bwz+cz^2&(\eqnumber)\cr
b^\prime=&\ b-2(awy-byz+cxz)&(\eqnumber)\cr
c^\prime=&\ ay^2-byx+cx^2.&(\eqnumber)\cr}
$$
We conclude that the greatest common divisor of $a$ and 2 is invariant under the 
action of $\overline{\Gamma}_0(2)$, hence $\overline{\Gamma}_0(2)$ acts on $\ca P^{\omega_2}_d$ and $\ca P^\omega_d$.
If $a$ is even, equation \eqn2  shows that the congruence class of $b$ modulo 4 
is also invariant under the action of $\overline{\Gamma}_0(2)$.
Hence for $k\in\{\pm 1\}$ the group $\overline{\Gamma}_0(2)$ acts on the subset
$\ca P^\omega_{d,k}=\{\tau\in\ca P_d: 2\mid a_\tau\hbox{ and }b_\tau\congr k\mod 4\}$ of 
$\ca P^\omega_d$.
We will use lemma  \thmn2 to prove the following three bijections: 
$$\displaylines{
\rlap{(\eqnumber)}\hfill \overline{\Gamma}_0(2)\backslash\ca P^{\omega_2}_d\ {\buildrel\sim\over\longrightarrow}\ 
\PSL_2(\bo Z)\backslash\ca P_d\hfill \cr
\noalign{\vskip2pt}
\rlap{(\eqnumber)}\hfill \overline{\Gamma}_0(2)\backslash\ca P^\omega_{d,k}\ {\buildrel\sim\over\longrightarrow}\ 
\PSL_2(\bo Z)\backslash\ca P_d\hfill  \llap{ for\ $k\in\{\pm 1\}$.\hskip1.5cm}\cr
}$$
In particular $\overline{\Gamma}_0(2)\backslash\ca P^{\omega_2}_d$ has cardinality $h$ and, because 
$\ca P^\omega_d=\ca P^\omega_{d,-1}\cup\ca P^\omega_{d,1}$, the set $\overline{\Gamma}_0(2)\backslash\ca P^\omega_d$ has 
cardinality $2h$.  
The set $I=\{\ma 1001 ,\ma 0{-1}10 ,\ma 1011 \}$ is a complete set of left coset 
representatives of $\overline{\Gamma}_0(2)$ in $\PSL_2(\bo Z)$. 
As $d=b^2-4ac$ is congruent to 1 modulo 8, we find for $\tau\in\ca P_d$ that $2|a_\tau c_\tau$ and 
$2\nmid b_\tau$.
In particular we have $\ca P^{\omega_2}_d=\{\tau\in\ca P_d:2\nmid a_\tau, 2\mid c_\tau\}$.
An easy calculation using \eqn5 and \eqn3 shows the equalities 
$\ma 0{-1}10 \ca P^{\omega_2}_d=\{\tau\in\ca P_d:2\mid a_\tau, 2\nmid c_\tau\}$
and 
$\ma 1011 \ca P^{\omega_2}_d=\{\tau\in\ca P_d:2\mid a_\tau, 2\mid c_\tau\}$.
Hence $\ca P_d$ is the disjoint union of $\{M\ca P^{\omega_2}_d:M\in I\}$ and \eqn2 follows from 
lemma  \thmn2.
In a similar way, the bijections \eqn1 follow from the descriptions
$$
\eqalign{
\textstyle\ca P^\omega_{d,k}= &\{\tau\in\ca P_d:2\mid a_\tau\hbox{ and } b_\tau\congr k\mod 
4 \},\cr
\textstyle\ma 0{-1}10 \ca P^\omega_{d,k}= &\{\tau\in\ca P_d:2\mid c_\tau\hbox{ and } 
b_\tau\congr -k\mod 4 \},\cr
\textstyle\ma 1011 \ca P^\omega_{d,k}= &\left\{\tau\in\ca P_d:\vcenter{
\hbox{$(2\nmid a_\tau\hbox{ and } b_\tau\congr k\mod4)$ or }
\hbox {$( 2\nmid c_\tau\hbox{ and } b_\tau\congr -k\mod 4)$}
}\right\}.\cr
}$$
As $d$ is congruent to 1 modulo 8 the Hilbert class field of $K$ coincides with the ray class field of conductor 2 of $K$.
Using \eqn{15} we find that each of the numbers $f({-1+\sqrt{d}\over 2})$ with $f\in\{\omega,\omega_1,\omega_2\} $  is an algebraic integer and generates the same field over $K$ as $j({-1+\sqrt{d}\over2})$.
In particular they are of degree $h$ over $K$.

It follows from the $q$-expansion of $\omega_2$  that $\omega_2({-1+\sqrt{d}\over2}) $ is real, hence of degree $h$ over $\bo Q$. 
Therefore its conjugates over $K$ coincide with those over $\bo Q$ and we can use \eqn{13} and \eqn{12}  to compute them.
Let $[ a,-b,c]$ be a quadratic form of discriminant $d\congr 1\mod 8$ with $a$ odd.
We find 
$$\textstyle{
\omega_2({-1+\sqrt{d}\over 2})^{[a,-b,c]}=\omega_2^{\ma 1{(b-1)/2}01}({-b+\sqrt{d}\over 
2a})=\omega_2({-b+\sqrt{d}\over 2a}),
}$$
where the last equality follows from the fact that $\omega_2$ is invariant under 
$\overline{\Gamma}_0(2)$.
As  ${-b+\sqrt{d}\over2a}\in\ca P^{\omega_2}_d$ and $\overline{\Gamma}_0(2)\backslash\ca P^{\omega_2}_d$ has cardinality $h$, we conclude that $\{\omega_2(\tau): 
\tau\in\overline{\Gamma}_0(2)\back\ca P_d^{\omega_2}\}$ is a complete set of conjugates of $\omega_2({-
1+\sqrt{d}\over 2})$ over $\bo Q $.

The conjugates of $\omega({-1+\sqrt{d}\over 2})$ and $\omega_1({-1+\sqrt{d}\over 2})$ over $K$
are computed similarly using the transformations  \eqn6.
For $[a,-b,c] \in\ca Q_d$ with $a$ odd we find 
$$\eqalign{
\omega(\textstyle{-1+\sqrt{d}\over 2})^{[a,-b,c]}&=\cases
\omega_2({b-2a+\sqrt{d}\over 2(c-b+a)})&\hbox{ if } b\congr 1\mod 4\cr
\omega_2({b+\sqrt{d}\over 2c})&\hbox{ if } b\congr 3\mod 4\cr
\endcases\cr
\noalign{\vskip2pt}
\omega_1(\textstyle{-1+\sqrt{d}\over 2})^{[a,-b,c]}&=\cases
\omega_2({b+\sqrt{d}\over 2c})&\hbox{ if } b\congr 1\mod 4\cr
\omega_2({b-2a+\sqrt{d}\over 2(c-b+a)})&\hbox{ if } b\congr 3\mod 4.\cr
\endcases\cr
}$$
Hence $\{\omega_2(\tau): \tau\in\overline{\Gamma}_0(2)\back\ca P_{d,1}^\omega\}$ is a complete set of 
conjugates of $\omega({-1+\sqrt{d}\over 2})$ over $K$ and 
$\{\omega_2(\tau): \tau\in\overline{\Gamma}_0(2)\back\ca P_{d,-1}^\omega\}$ is a complete set of conjugates 
of $\omega_1({-1+\sqrt{d}\over 2})$ over $K$.
Recall that $\omega_2 $ maps $\overline{\Gamma}_0(2)\back\bH $ bijectively to $\bo P ^1(\bo C)$.
In order words for $\tau,\tau^\prime\in\H$ the equation $\omega_2(\tau)=\omega_2(\tau^\prime)$ holds if and only 
if $\tau=\gamma\tau^\prime$ for some $\gamma\in\overline{\Gamma}_0(2)$. 
Therefore $\omega({-1+\sqrt{d}\over 2})$ and $\omega_1({-1+\sqrt{d}\over 2})$ are not 
conjugate over $K$. From the  $q$-expansion  we see that $\omega({-
1+\sqrt{d}\over 2})$ and $\omega_1({-1+\sqrt{d}\over 2})$ are complex conjugates and 
hence conjugate over $\bo Q$. 
We conclude that $\{\omega_2(\tau): \tau\in\overline{\Gamma}_0(2)\back\ca P_d^\omega\}$ is a complete set 
of $2h$ conjugates over $\bo Q$ of both $\omega({-1+\sqrt{d}\over 2})$ and  $\omega_1({-
1+\sqrt{d}\over 2})$.\hfill$\qed$
\bigskip \noindent
Similar to \eqn7 we find for $f\in\{\omega,\omega_1,\omega_2\} $ 
the following corollary which we will use in the next section:
$$\eqalign{
{\textstyle\log \N\big ({\textstyle f({-1+\sqrt{d_1}\over2}),f({-1+\sqrt{d_2}\over2})}\big)}&
\hskip 0.1 cm
=
\hskip -0.9 cm
\sum_{i=1,2,\tau_i\in\overline{\Gamma}_0(2)\backslash\ca P_{d_i}^{f}}
\hskip -0.7cm
\log |\omega_2(\tau_1) -\omega_2(\tau_2) |\cr
&\hskip -2.4cm
=
\hskip -0.9 cm
\sum_{i=1,2,\tau_i\in\overline{\Gamma}_0(2)\backslash\ca P_{d_i}^{f}}
\hskip -0.7 cm
 {\textstyle \log |{2^{12}\over\omega_2(\tau_1)} -{2^{12}\over\omega_2(\tau_2)} |}-
\cases
12h_1h_2\log 2 &\text{\rm  if }f =\omega_2\cr
0&\text{\rm  otherwise.}\cr
\endcases 
}
\leqno(\eqnumber)
$$
To obtain the last line we used that the absolute value of the norm of $\omega({-1+\sqrt{d_i}\over2})$ and $\omega_2({-1+\sqrt{d_i}\over2})$ 
is equal to $2 ^{12h_i} $ and $1$ respectively [\LA, chapter 12 \S 2].

\head 3. A limit formula for differences of hauptmodul values
\endhead
\def\P#1{#1}
\noindent\rm To describe the results of this section, we concentrate for the moment on the function $\gamma_2$; similar remarks hold for the other functions.
The main result for  $\gamma_2$  is theorem 13 below.
With $\alpha_1,\alpha_2$ as in theorem 4,
it describes the absolute norm $\N\big(\gamma_2(\alpha_1),\gamma_2(\alpha_2)\big)$
in terms of a limit value of a certain meromorphic function.
This meromorphic function is given as an infinite sum over {\it integers}.
The proof of this result is in two steps.  
The first step is a general function theoretic result, independent of the results in the previous section.
For all $ z_1, z_2\in\bo H $ such that  $ z_1\not\in\Gamma^3 z_2$ we express
$\log |\gamma_2( z_1) -\gamma_2( z_2) |$ as the constant term of the Laurent expansion at 1 of a meromorphic function.
Here the meromorphic function is defined as an infinite sum over {\it matrices} $\gamma\in\Gamma^3$.
The second step is more algebraic.
We substitute $z_i =\alpha_i \ (i = 1,2)$ and sum over the conjugates of 
$\gamma_2(\alpha_1) $ and $\gamma_2(\alpha_2)$.
The description of these conjugates in section 2 enables us to transform the infinite sum over matrices
into an `easier' infinite sum over integers.
\bigskip\noindent
The group of complex analytic automorphisms of $\bo H$ is isomorphic to $\PSL_2(\bo R)$.
Its elements, which we represent by matrices, act on $\bo H$ by linear fractional transformations.  
Although we are mainly interested in the subgroups $\P{\Gamma ^ 3} $ and $\overline{\Gamma}_0(2) $
of $\PSL_2(\bo R)$,
our setup will be more general.
Let $\Gamma\subset\PSL_2(\bo R)$ be a discrete subgroup and let
$\ca C_\Gamma$ be the set of cusps of $\Gamma$.
Recall that the cusps of $\Gamma$ are those $z\in\bo R\cup\{\infty\}$ that are
the unique fixed point of some $\gamma\in\Gamma$.
The action of $\Gamma$ on  $\bo H$ extends to an action on $\bH=\bo H\cup\ca C_\Gamma$
and we endow  $\Gamma\back\bH$ with the usual Riemann surface structure [\MI, \S 1.8].
The  $\PSL_2(\bo R)$-invariant measure ${dxdy\over y^2}$ on $\bo H$ gives rise to a measure $\mu$ on
 $\Gamma\back\bH$. 
As an example we consider the group $\Gamma =\PSL_2(\bo Z)$.
The cusps are the elements of $\bo Q\cup\{\infty\}$ and  $\PSL_2(\bo Z)$ acts transitively on this set. Moreover $\PSL_2(\bo Z)\back\bH$ has genus zero and $\mu$-measure ${\pi\over3}$ [\MI, \S 4.1].
 
We make the following three assumptions on the group $\Gamma $.
First of all we assume that  $\Gamma\back\bH$ is compact.
By a theorem of Siegel [\MI, page 32] this is equivalent to the assumption that  $\Gamma\back\bH$ 
has finite volume.
This implies in particular that $\Gamma$ is finitely generated [\IW, page 41] and that the set $\ca C_\Gamma$
has finitely many $\Gamma$-orbits [\MI, page 27].
Furthermore we want $\Gamma$ to have at least one cusp.  To fix notation we take
$\infty\in\ca C_\Gamma$.
Finally we suppose that the Riemann surface $\Gamma\back\bH$ has genus zero.
There are infinitely many conjugacy classes of groups satisfying these three conditions.
In this paper we will focus on the groups $\PSL_2(\bo Z)$,
$\P {\Gamma^3}$ and $\overline{\Gamma}_0(2)$.
In the previous paragraph we noted that $\PSL_2(\bo Z)$ does indeed satisfies the three conditions.
The other two groups are of finite index in $\PSL_2(\bo Z) $ and therefore have finite co-volume
and $\infty$ as a cusp.
The fact that $\Gamma^3\back\bH$ and $\overline{\Gamma}_0(2)\back\bH$ are of genus zero was already established in section 2.

As $\Gamma\back\bH$ has genus zero by assumption, its function field is generated by one element.
To single out a special kind of generator we first recall that the stabilizer
 $\Gamma_\infty=\{\gamma\in\Gamma:\gamma\infty =\infty\} $ 
of the cusp $\infty$ is an infinite cylic group generated by $\ma 1h01$ for some positive $h$.
The scaling transformation $\sigma_\infty =\ma {\sqrt{h}}00{1\over\sqrt{h}}$ fixes $\infty$ and  $\sigma_\infty^{-1}\Gamma_\infty\sigma_\infty$ is generated by the translation $z\mapsto z+1$.
Hence for any generator~$f$ of the above function field
the meromorphic function $f(\sigma_\infty z)$ has a Laurent expansion in $q=e^{2\pi i z}$. 
Such a generator $f$ is called a {\it normalised principal modulus} (or {\it normalised hauptmodul}) for $\Gamma$ if this expansion has the form
$$
f(\sigma_\infty z) =q^{-1} +\sum_{n\geq 0} a_nq^n,\leqno(\eqnumber)
$$
with $a_n\in\bo C $.  
In particular a normalised principal modulus has its pole at $\infty$ and is unique up to an additive constant.  
We will consider $f$ both as a function on the Riemann surface $\Gamma\back\bH$ and 
as a $\Gamma$-invariant function on $\H$.
For example $\gamma_2$ is a normalised principal modulus for the group $\P{\Gamma^{3}} $  
and $-2^{12}/\omega_2 $ is a normalised principal modulus for $\overline{\Gamma}_0(2)$. 

Theorem  \thmn{-1}  below states an equality between 
$\log | f( z_1) - f( z_2) |$ and the limit value of a certain meromorphic function.
This equality will follow from the following characterization of the function $\log | f( z_1) - f( z_2) |$.
\bigskip\noindent\bf Proposition\,\thmnumber\sl 
Let $\Gamma\subset\PSL_2(\bo R) $ be as above and let $f$ be a normalised principal modulus for $\Gamma $.
The function 
$$h( z_1, z_2) =\log | f( z_1) - f( z_2) |^{2}$$
is the unique symmetric function on
$\bo H\times\bo H -\{( z_1, z_2):  z_1\not\in\Gamma z_2 \} $ 
that for fixed $ z_2\in\bo H$ satisfies 
\item{1.} $h(\gamma z_1, z_2)=h( z_1, z_2)$  for all $\gamma\in\Gamma$; 
\item{2.} $({\partial^2\over\partial x_1^2}+{\partial^2\over\partial y_1^2})h( z_1, z_2) = 0$, where $x_1 =\Re( z_1)$
and $y_1 =\Im( z_1)$; 
\item{3.} $h( z_1, z_2) = e_{ z_2}\log |  z_1- z_2 |^{2} +O(1)$ for $z_1\rightarrow z_2$,  
where $e_{ z_2}$ denotes the order of the stabilizer of $ z_2$ in $\Gamma$;
\item{4.} $h(\sigma_\infty z_1, z_2) = 4\pi\Im( z_1) +o(1)$ for $ z_1\rightarrow\infty$; 
\item{5.}  $h( z_1, z_2) = O(1)$ for $ z_1\rightarrow a$ with $a\in\ca C_\Gamma $ not $\Gamma$-equivalent to $\infty$.
\bigskip \noindent \bf Proof. \rm
The first two properties are obvious.
To prove property 3 we use the following equality:
$$
f( z_1) -f( z_2) =( z_1- z_2)^{e}\cdot { f( z_1) -f( z_2) \over ( z_1- z_2)^{e} }.
$$
with $e=e_{z_2}$.
For fixed $ z_2\in\bo H$ the function $f( z_1) -f( z_2) $ has a zero of order $e_{ z_2}$ in $ z_1 = z_2$.
Therefore, the limit for $ z_1\rightarrow z_2 $ of the second factor on the right converges
to a non-zero value and the third property follows.
Property~4 follows by using the  Fourier expansion at infinity \eqn1 of the function  $f$. 
For property~5 we use that $f$ is a bijection between $\Gamma\back\bH$ and the projective line over  $\bo C$.
If $ z_2\in\bo H$ is fixed, we find that $f(z_1)\not= f( z_2)$ for $ z_1$ close to a cusp $a$.
By definition $f$ has its pole at  $\infty$ , hence for every cusp $a $ not $\Gamma$-equivalent to $\infty$ 
the value $| f( z_1) -f( z_2) | $ is nonzero and bounded from above for 
$z_1\rightarrow a$ and property 5 follows. 

To prove uniqueness, we let  $k( z_1, z_2)$ be the difference of two functions on 
$\bo H\times\bo H-\{( z_1, z_2):  z_1\not\in\Gamma z_2\} $ that satisfy the above five properties.
Fix $ z_2\in\bo H $ and consider $k( z_1, z_2)$ as a function on  
$\bo H-\{ z_1\not\in\Gamma z_2\} $.
On this domain the function is  $\Gamma$-invariant and harmonic by the first two properties.
According to~3,4 and 5 the function $k( z_1, z_2) $ is bounded if $ z_1$ approaches
an element of $\Gamma  z_2$  or a cusp of $\Gamma$.
Hence $k( z_1, z_2)$ extends to a harmonic function on the Riemann surface $\Gamma\back\bH$.
As $\Gamma\back\bH $ is compact by assumption we find that  $k( z_1, z_2)$ is constant. 
Using the fourth property we conclude that the function $k$ is identically zero, which proves the 
uniqueness of $h$.
\hfill $\qed$

\medskip\noindent
We will now `build' a function from scratch which satisfies the five properties above.
To motivate the definitions below, we follow the arguments in [\GZt, II \S 2].
First we make a bi-$\Gamma$-invariant function $G( z_1, z_2) $ on $\bo H\times\bo H $.
Let $g( z_1, z_2)$ be a function on $\bo H\times\bo H $ 
satisfying $g(\gamma z_1,\gamma z_2) =g( z_1, z_2) $ for all $\gamma\in\PSL(\bo R)$.
This is equivalent with $g$ being a function of the hyperbolic distance $d(z_1,z_2)$
between $ z_1$ and  $ z_2$. In particular $g$ is symmetric.  
Ignoring convergence, the function $G( z_1, z_2) =\sum_{\gamma\in\Gamma}g( z_1,\gamma z_2) $ is then clearly bi-$\Gamma$-invariant.     
If we want~$G$ to satisfy properties~2 and~3 it seems natural to select a $g$ that is harmonic in both variables
and satisfies $g( z_1, z_2) =\log | z_1- z_2 |^{2}+ O(1) $ for $ z_1\rightarrow z_2$. 
Unfortunately this leads to convergence problems in the definition of $G$.
To resolve this difficulty we weaken the condition that $g$ should be harmonic.
For a complex variable  $z =x+iy\in\bo H$,  let $\Delta=\Delta_z =y^{2}({\partial^2\over\partial x^2} + {\partial^2\over\partial y^2} )$ denote the hyperbolic Laplacian in  $x,y$.
Recall that $\Delta$ is an invariant operator:  $\Delta (f(\gamma z)) =(\Delta f)(\gamma z)$ for all sufficiently smooth functions $f$ on $\bo H$ and $\gamma\in\PSL(\bo R) $.   
We choose $g$ to be `almost' harmonic, in the sense that $g$ satisfies 
$\Delta g=\epsilon g $  for some small  $\epsilon>0$, 
and with a logarithmic singularity along the diagonal as above.
Note that as $g$ is symmetric, we can take either $\Delta=\Delta_{ z_1}$ or  $\Delta=\Delta_{ z_2}$.  
Because $g( z_1, z_2) $ is a function of the hyperbolic distance $d( z_1, z_2) $, we can write
$g( z_1, z_2) =Q\big (1+{|z_1-z_2|^2\over 2\Im(z_1)\Im(z_2)}\big)$ for some function  $Q$,
the argument of $Q$ being $\cosh\big (d(z_1,z_2)\big)$.
The partial differential equation $\Delta g =\epsilon g$ translates into the ordinary differential 
equation 
$$\textstyle
\big ( (1-t ^2){d ^2\over dt ^2} -2t{d\over dt} + \epsilon\big)Q(t) =0
$$
for the function $Q$.
This is the Legendre differential 
equation of index $s-1 $, where $\epsilon=s(s-1)$ with  $s > 1 $.
Up to a scalar it has a unique solution 
that is small at infinity, the Legendre function of the second kind  $Q_{s-1}(t) $.
This function,  real analytic in $ t\in\bo R_{>1}$ and holomorphic in $s\in\bo H_1=\{s\in\bo C:\Re(s) > 1\}$, is given by 
$$\textstyle
Q_{s-1}(t)={\Gamma(s)^2\over 2\Gamma(2s)} \big ({2\over 1+t}\big)^s  F(s,s ;2s ;{2\over 1+t}), \leqno(\eqnumber)
$$  
where $F(a,b;c;z)$ is Gauss's hypergeometric function and $\Gamma(s)$ is the gamma function
(in the sequel it will be clear from the context whether $\Gamma$ denotes a group
or the gamma function).
The Legendre function of the second kind has the following asymptotic behaviour [\IW, lemma 1.7]
$$
\displaylines{
\rlap{(\eqnumber)}\hfill Q_{s-1}(t) = -\half\log(t-1) + O(1)\hskip0.3cm\text{\rm  for }t\downarrow 1,\hfill \cr
\rlap{(\eqnumber)}\hfill Q_{s-1}(t) = O(t^{-s})\hskip0.3cm \text{\rm for }t\rightarrow\infty.\hfill \cr
} $$
The above discussion leads to the following definitions.
For  $(z_1, z_2,s)\in\bo H\times\bo H \times \bo H_1$ with $z_1\not=z_2 $ the function
$$\leqalignno{\textstyle
g_s(z_1,z_2)=&\textstyle -2Q_{s-1}\big(1+{|z_1-z_2|^2\over 2\Im(z_1)\Im(z_2)}\big) &(\eqnumber)\cr
}$$
satisfies
$$\displaylines {
\rlap{(\eqnumber)}\hfill  g_s(\gamma z_1,\gamma z_2) =g_s( z_1, z_2)\hskip0.3cm \text{\rm for all }
\gamma\in\PSL(\bo R)\hfill\cr
\rlap{(\eqnumber)}\hfill \Delta_{z_i} g_s( z_1, z_2) =s(s-1)g_s( z_1, z_2)\hskip0.3cm (i = 1, 2) .\hfill\cr 
}$$
For $(z_1,z_2,s)\in\bo H\times\bo H \times\bo H_1 $ with $ z_1\not\in\Gamma z_2$ we define the {\sl automorphic Green function}, or {\sl resolvent kernel function},  for $\Gamma$ by
$$
G_{\Gamma,s}(z_1,z_2)=\sum_{\gamma\in\Gamma}g_s(z_1,\gamma z_2),\leqno(\eqnumber)
$$
which is $4\pi $ times the function studied in [\HEj] and [\IW].
Using \eqn5 one shows that this sum
is uniformly and absolutely convergent on compact subsets of its domain of definition~[\HEj, page 31].
Consequently $G_{\Gamma,s}(z_1,z_2)$ is holomorphic as a function of $s\in\bo H_1 $ and, because of \eqn6 and \eqn3,  satisfies property 1 and 3 of proposition~\thmn1:
$$\displaylines {
\ \llap(\eqnumber)\hfill G_{\Gamma,s}(\gamma  z_1, z_2) =G_{\Gamma,s}( z_1, z_2)
\hskip0.3cm\text{\rm for } \gamma\in\Gamma\hfill\cr
\ \llap(\eqnumber)\hfill G_{\Gamma,s}( z_1, z_2) = e_{ z_2}\log |  z_1- z_2 |^{2} +O(1)
\hskip0.3cm\text{\rm for fixed } z_2\in\bo H\hskip0.15cm\text{\rm and } z_1\rightarrow z_2\hfill\cr
}$$
\noindent
Using that  \eqn7 is also valid for derivatives of every order of $Q_{s-1}(t)$ one can show that $G_{\Gamma,s}(z_1,z_2)$
is infinitely differentiable as a function of the real variables $x_1,y_1,x_2 $ and $y_2$,
provided that $z_1$ and $z_2$ are not equivalent modulo $\Gamma$. 
Furthermore, by equation \eqn4 and the fact that the Laplacian is an invariant operator we have
$$
\Delta_{z_i}G_{\Gamma,s} ( z_1, z_2) =s(s-1)G_{\Gamma,s} ( z_1, z_2)\ \ \ (i=1,2).\leqno(\eqnumber)
$$
\noindent
By taking the limit $s\rightarrow 1$ of $G_{\Gamma,s}( z_1, z_2) $ one might hope to 
obtain a harmonic function. Unfortunately this limit does not exist.
The function $G_{\Gamma,s}(z_1,z_2)$ does have a meromorphic continuation to the entire $s$-plane
which satisfies \eqn3, \eqn2 and \eqn1.   
At $s = 1 $ it has a simple pole with residue $-4\pi\mu(\Gamma\back\bo H)^{-1} $ independent of $z_1$
and~$z_2$ [\HEj, chapter 8.6; \IW, chapter 7.4].
Therefore we first subtract the `singular part of $G_{\Gamma,s}( z_1, z_2)$ at  $s = 1 $'
and then take the limit  $s\rightarrow 1 $.
As we will see below this results in a function which is not only harmonic but satisfies all five properties of proposition~\thmn1.

To obtain this `singular part'  we fix $ (z_2,s)\in\bo H\times\bo H_1$,  and consider
$G_{\Gamma,s}( z_1, z_2)$ as an infinitely differentiable function of
$x_1$ and $y_1$ with $y_1 > \text{\rm sup}\{\Im(\gamma  z_2):\gamma\in\Gamma \}$.
The function $G_{\Gamma,s}(\sigma_\infty   z_1, z_2) $  
is invariant under $x_1\mapsto  x_1+1 $, hence it has a Fourier expansion in $e^{2\pi i x_1}$.
One can show [\HEj, (3.3) on page 274] that its zeroth Fourier coefficient is equal to  ${4\pi\over 1-2s}E_\Gamma(\sigma_\infty  z_2,s) y_1^{1-s}$, with
$$
E_\Gamma( z,s)=\sum_{\gamma\in\Gamma_\infty\backslash\Gamma}\Im(\sigma_\infty^{-1}\gamma 
 z)^s\leqno(\eqnumber)
$$
the Eisenstein series for $\Gamma$ and the cusp $\infty$. 
For the modular group $\PSL_2(\bo Z)$ this function can be written in the more familiar form
$$
E_{\PSL_2(\bo Z)}(z,s) =\half\sumprime_{c,d\in\bo Z\atop(c,d) =1} {y ^s\over | cz +d |^{2s}}.
$$
The sum  \eqn1 converges absolutely and uniformly on compact subsets of $\bo H\times\bo H_1$, 
is real analytic as a function of $z\in\bo H $ and is holomorphic as a function of  $s\in\bo H_1$.
In fact one can prove that  $E_\Gamma(z,s) $ admits a meromorphic continuation to the whole  $s$-plane
 with a simple pole at $s = 1 $ with residue $\mu(\Gamma\back\bo H)^{-1} $
independent of~$z$~[\IW, chapter~6].
At regular points $s\in\bo C $ the Eisenstein series satisfies the following properties:
$$\displaylines {
\ \llap(\eqnumber)\hfill E_\Gamma(\gamma z,s)=E_\Gamma(z,s)\ \ \ \text{\rm  for all }\gamma\in\Gamma\hfill \cr
\ \llap(\eqnumber)\hfill \Delta E_\Gamma(z,s)=s(s-1)E_\Gamma(z,s).\hfill \cr
} $$
Apart from the zeroth one the Fourier coefficients of $G_{\Gamma,s}(\sigma_\infty z_1,z_2)$ are all regular at $s =1$ [\HEj, chapter 8.6].
For fixed $ z_2\in\bo H$ and $y_1\rightarrow\infty$ we have
$$
\lim_{s\rightarrow1}\textstyle\left[ G_{\Gamma,s}(\sigma_\infty  z_1, z_2)- {4\pi\over 1-2s}E_\Gamma(  z_2,s)y_1^{1-s}\right] =O(e^{-2\pi y_1}),\leqno(\eqnumber)$$
which will be used in the proof of theorem \thmn{-1} below.

If we consider $E_\Gamma(\sigma_\infty  z,s) $ as function of $ z =x+iy\in\bo H$, 
it is invariant under the translation $x\mapsto  x+1 $.
Its zeroth Fourier coefficient is equal to $y^s+\phi_\Gamma(s)y ^{1-s}$ [\IW, page 66] with  
$$
\phi_\Gamma(s)= {\textstyle\sqrt\pi{\Gamma(s-{1\over 2})\over \Gamma(s)}} \sum_{c>0}c^{-2s}
\textstyle\#\left\{d\mod c:\ma **cd \in\sigma_\infty^{-1}\Gamma\sigma_\infty\right\}\leqno(\eqnumber)
$$
for $s\in\bo H_1$.
For the group $\PSL_2(\bo Z)$ the above sum is over positive integers $c$ and $\sigma_\infty$ is the identity transformation.
For $c\in\bo Z_{> 0} $  the cardinality of the set in \eqn1 is equal to the number of positive integers less than and coprime to $c$.  We obtain
$$
\phi_{\PSL_2(\bo Z)}(s) =\sqrt{\pi}{\Gamma(s-{1\over 2})\over \Gamma(s)}{\zeta(2s-1)\over\zeta(2s)}
$$
with $\zeta(s)$ the Riemann zeta function.

The function $\phi_\Gamma(s)$ has a meromorphic continuation to $\bo C $ with the same behaviour
at $s =1 $ as $E_\Gamma(\sigma_\infty  z,s)$.
More precisely, for $y\rightarrow\infty $ 
$$
\lim_{s\rightarrow1}\textstyle\left[
E_\Gamma(\sigma_\infty z,s) - \big(y^s+ \phi_\Gamma(s)y^{1-s} \big)\right] =
O(e^{-2\pi y}).\leqno(\eqnumber)
$$
The following result is a slight generalization of proposition 5.1 in [\GZe].
For a normalised principal modulus $f$ and $z_1,z_2\in\bo H $ it gives a formula for the archimedean local height of the degree zero divisors $(f(z_1)) -(\infty) $ and $(f(z_2)) -(\infty)$ on $\bo P ^1$ [\GR; \GZd,~IV~\S 3; \GZt,~II~\S 2].
\bigskip\noindent\bf Theorem\,\thmnumber\sl
Let $\Gamma\subset\PSL_2(\bo R)$ be a discrete subgroup for which $\infty$ is a cusp
and for which the Riemann surface $\Gamma\back\bH$ is compact of genus zero.
Furthermore let $f$ be a normalised principal modulus for $\Gamma$.
For $ z_1, z_2\in\bo H$ with $ z_1\not\in\Gamma z_2$ we have the following equality:
$$
\log|f( z_1)-f( z_2)|^2=\lim_{s\rightarrow 1}
\Big[ G_{\Gamma,s}(z_1,z_2)-\textstyle{{4\pi\over 1-2s}}
\big(E_\Gamma( z_1,s)+E_\Gamma( z_2,s)-\phi_\Gamma(s)\big)\Big]
$$
\bigskip \noindent\bf Proof.
\rm It follows from the above discussion that the function between square brackets on 
the right has at most a simple pole at $s = 1 $.
As the residues add up to zero, we conclude that the limit exists.

Let $h( z_1, z_2)$ denote the function on the right of the equality sign above.
We show that this symmetric function satisfies the five properties of proposition  \thmn2, from which the theorem follows.
Fix $ z_2\in\bo H $ and consider $h(z_1,z_2)$ as a function of~$z_1$.
The fact that $h( z_1, z_2)$ is $\Gamma$-invariant is obvious from \eqn9 and \eqn5.
For $s$ in some punctured neighbourhood of  $s = 1 $ 
$$
\textstyle{G_{\Gamma,s}(z_1,z_2)-{4\pi\over 1-2s}E_\Gamma( z_1,s)}
$$
is a $\Delta_{z_1}$-eigenfunction with eigenvalue  $s(s-1)$, because of \eqn7 and \eqn4.
In particular its limit  $s\rightarrow 1$ is harmonic. 
The remaining terms of  $h$ are independent of $ z_1$, so we conclude that $h$ is harmonic. 
Using \eqn8, we find
$$h( z_1, z_2) = e_{ z_2}\log |  z_1- z_2 |^{2} +O(1)\text{\rm \ \ for\ }z_1\rightarrow z_2,$$
with $e_{ z_2}$ the order of the stabilizer of $ z_2$ in $\Gamma$.

To study $h(z_1,z_2)$ at the cusp $\infty$ we use the equations \eqn3 and \eqn1 (with $z=z_1$), which are equivalent to 
$$\displaylines{
\lim_{s\rightarrow 1}[G_{\Gamma,s}(\sigma_\infty z_1,z_2)-\textstyle{4\pi\over 1-2s}E_\Gamma( z_2,s)]=
\textstyle{4\pi\over\mu(\Gamma\back\bo H)}\log y_1+ O(e^{-2\pi y_1})\cr
\noalign{\noindent and}
\lim_{s\rightarrow 1}\textstyle{4\pi\over 1-2s}[E_\Gamma(\sigma_\infty z_1,s) -\phi_\Gamma(s)]=
-4\pi y_1+\textstyle{4\pi\over \mu(\Gamma\back\bo H)}\log y_1+O(e^{-2\pi y_1}),\cr
} $$
respectively.
By subtracting the last equation from the first we conclude that $h( z_1, z_2)$ satisfies the fourth property of proposition  \thmn2. 

In order to examine $h( z_1, z_2) $ for $ z_1$ near a finite cusp, we need to generalize the estimates
\eqn3 and \eqn1 to all cusps. 
Fix a cusp $a\in\ca C_\Gamma$.
Choose a scaling transformation $\sigma_a\in\PSL_2(\bo R)$ such that
$$\textstyle
\sigma_a\infty = a\hskip0.5cm\text{\rm and }\hskip0.5cm \sigma_a^{-1}\Gamma_a\sigma_a = \langle\ma 1101\rangle. 
$$
This transformation is unique up to composition on the right with a translation.
We need to show that the limit of $h(\sigma_a z_1,z_2)$ for $y_1\rightarrow\infty$ exists if $a$ is not $\Gamma$-equivalent to $\infty$.  
By comparing $G_{\Gamma, s} ( \sigma_a z_1 , z_2)$ with its zeroth Fourier coefficient we find for $y_1\rightarrow\infty $
$$
\lim_{s\rightarrow1}\textstyle\left[ G_{\Gamma,s}(\sigma_a z_1, z_2)- {4\pi\over 1-2s}E_{\Gamma,a} (  z_2,s)y_1^{1-s}\right] =O(e^{-2\pi y_1}).\leqno(\eqnumber)
$$
The Eisenstein series $E_{\Gamma,a}(z,s)$ for $\Gamma$ and the cusp $a$ is  for $s\in\bo H_1$ defined by
$$
E_{\Gamma,a} (z,s)=\sum_{\gamma\in\Gamma_\infty\backslash\Gamma}\Im(\sigma_a^{-1}\gamma 
z)^s.
$$
As in the case  $a =\infty $ (see  \eqn7) this function has a meromorphic continuation to the whole $s$-plane 
with a simple pole at  $s = 1 $ with residue  $\mu(\Gamma\back\bo H)^{-1}$. 

The Fourier expansion of $E_\Gamma(\sigma_a z_1,s) =E_{\Gamma,\infty}(\sigma_a z_1,s) $ 
leads for $y_1\rightarrow\infty$ to
$$
\lim_{s\rightarrow1}\textstyle\left[
E_\Gamma(\sigma_a z_1,s) - \big(\delta_{a\infty}y_1^s+ \phi_{\Gamma,a}(s) y_1^{1-s} \big)\right] =
O(e^{-2\pi y_1}).\leqno(\eqnumber)
$$
%
Here  $\delta_{a\infty}$ is 1 if $a$ and $\infty$ are $\Gamma$-equivalent, and 0 otherwise.
Furthermore $\phi_{\Gamma, a}(s)$ is a meromorphic
function in $s$ with a simple pole at $s = 1 $ with the same residue as the Eisenstein series. 
Note that $\phi_{\Gamma,\infty}(s) =\phi_\Gamma(s)$ is defined by  \eqn4.
Using  \eqn2 and \eqn1  and the fact that both $E_{\Gamma,a}(z,s) -E_\Gamma(z,s) $ and $\phi_{\Gamma,a}(s) -\phi_\Gamma(s)$ are bounded for $s\rightarrow1$,
 we can compute the limit of $h(\sigma_a z_1, z_2)$ for $ y_1\rightarrow \infty$.
This limit is finite if $a $ is not  $\Gamma$-equivalent to $\infty $  because of the  $\delta_{a\infty}$ in~\eqn1.
\hfill $\qed $
\bigskip\rm\noindent 
Let $d_1$ and $d_2$ be negative fundamental discriminants that are relatively prime and both congruent to~2 modulo~3.
For $i =1,2 $  let $\alpha_i\in\H$ be of discriminant~$d_i$ and normalised as in theorem~\thmn7.
We will apply the previous theorem  for the group~$\P{\Gamma ^3}$ and its normalised principal modulus $\gamma_2$ 
together with  \eqn{27}  to obtain a formula for 
$\log \N\big(\gamma_2(\alpha_1),\gamma_2(\alpha_2)\big) $.
Therefore we sum each of the four terms on the right hand side of the equation in theorem~\thmn1 
over $z_i\in\Gamma ^3\back\ca P_{d_i}^{\gamma_2}$, $i =1,2$.  

We start with the automorphic Green function.
Replace the summation variable $\gamma\in\P{\Gamma^3} $ in the definition \eqn{12}  of  $G_{\P{\Gamma^3},s}( z_1, z_2)$
by $\kappa_1^{-1}\kappa_2$ with   $\kappa_1,\kappa_2\in\P{\Gamma^3} $.
The elements $\kappa_1,\kappa_2$  are well defined up to right multiplication by elements of the stabilisers
 $\P{\Gamma^{3}_{z_1}},\P{\Gamma^{3}_{z_2}} $ and up to simultaneous left multiplication by an element of $\P {\Gamma^{3}}$.  
Using that the stabilisers $\P{\Gamma^3_{z_i}} $ for $z_i\in\ca P_{d_i}^{\gamma_2}$ have order $w_i\over 2$ 
and the invariance~\eqn{14} of the function $g_s( z_1, z_2)$  we find similarly to [\GZe, page 209]
the following:
$${\textstyle{4\over w_1w_2}}\hskip-0.4cm\sum_{i=1,2,\tau_i\in\Gamma^3\backslash\ca P_{d_i}^{\gamma_2}}
G_{\P{\Gamma^3},s}(\tau_1,\tau_2)=
\hskip -0.2cm\sum_{(\tau_1,\tau_2)\in\Gamma^3\back\ca P_{d_1}^{\gamma_2} \times\ca P_{d_2}^{\gamma_2}}
g_s(\tau_1,\tau_2),\leqno(\eqnumber)
$$
where $\Gamma^3 $ acts diagonally on $\ca P_{d_1}^{\gamma_2} \times\ca P_{d_2}^{\gamma_2}$. 
By definition  \eqn{16}  we have for $\tau_i ={-b_i+\sqrt{d_i}\over 2a_i} $  the equality
$$\textstyle
g_s(\tau_1,\tau_2) = -2Q_{s-1}\big(\textstyle{B(\tau_1,\tau_2)\over \sqrt{D}}\big) 
$$
with $D =d_1d_2$ and where
$$
B(\tau_1,\tau_2) =2a_1c_2+2a_2c_1 -b_1b_2
$$
is an integer which is larger than $\sqrt{D}$ by  \eqn{16}.
Fix an integer $n >\sqrt{D}$.
Because of the absolute convergence of  \eqn1, the number of $(\tau_1,\tau_2)$ on the right hand side of
 \eqn1 for which the argument of $Q_{s-1}$ equals $n\over\sqrt{D} $  is finite.
We obtain
$$
{\textstyle{4\over w_1w_2}}\hskip-0.4cm\sum_{i=1,2,\tau_i\in\Gamma^3\backslash\ca P_{d_i}^{\gamma_2}}
G_{\P{\Gamma^3},s}(\tau_1,\tau_2)=
-2\hskip -0.1cm\sum_{n>\sqrt{D}}\rho^{\gamma_2}(n)Q_{s-1}(\textstyle{n\over\sqrt{D}})\leqno(\eqnumber)
$$
where $\rho^{\gamma_2}(n)$ is the cardinality of the set
$$
\ca S^{\gamma_2}_{d_1,d_2,n}=\Gamma^3\backslash\{(\tau_1, 
\tau_2)\in\ca P^{\gamma_2}_{d_1}\times\ca P^{\gamma_2}_{d_2}: B(\tau_1,\tau_2)=n\}.
$$
For integers  $n$ with $n^2\not\congr D \mod 36 $ this set is empty.
Namely for $\tau_i ={-b_i+\sqrt{d_i}\over 2a_i}$ we have
$$
B(\tau_1,\tau_2)^2-D =4(a_1c_2- a_2c_1)^2+4(a_1b_2-a_2b_1)( c_1b_2-c_2b_1).\leqno(\eqnumber)
$$
If $\tau_i\in\ca P_{d_i}^{\gamma_2}$ then $b_i$ is divisible by 3 and  $a_i\congr c_i\mod 3$, hence \eqn1 
is divisible by~36. 

To determine  $\rho^{\gamma_2}(n)$ for integers $n$ that satisfy $n^2\congr D \mod 36 $  we use the following proposition from [\GZd].
\bigskip\noindent\bf Proposition\,\thmnumber\sl
Let $N$ be a positive integer, let $d_1$ and $d_2$ be negative relatively prime 
integers such that $d_1$ and $d_2$ are squares modulo $4N$ and set $D =d_1d_2 $.
For positive integers $n$ define the set
$$
\ca S_{d_1,d_2,n} ^N =\overline {\Gamma}_0(N)\backslash\{(\tau_1, \tau_2)\in\ca P_{d_1} ^N\times
\ca P_{d_2} ^N:B(\tau_1,\tau_2)=n\}
$$
with $\ca P_{d_i} ^N =\{\tau\in\ca P_{d_i}:N\mid a_{\tau}\}$ for $i =1,2$.
The set $\ca S_{d_1,d_2, n} ^N $ has cardinality
$$
\#\ca S_{d_1,d_2, n} ^N =\cases
2^t\sum_{d|{n^2-D\over 4N}}\varepsilon(d)&\text{\rm  if }n^2\congr D\mod 4N;\cr
0&\text{\rm  otherwise,}\cr
\endcases$$
where $t$ is the number of primes dividing $N$ and $\varepsilon$ is defined as in the introduction.
\bigskip \noindent\bf Proof. \rm The case $N =1 $ is proved in proposition 6.1 of [\GZe] using class field theory.
The general case is proved using the theory of quaternion algebras on page 516 of [\GZd].
\hfill $\qed $
\bigskip\noindent\rm
We return to the situation before proposition \eqn1, i.e.  $d_1$ and $d_2 $ are coprime negative fundamental discriminants both congruent to 2 modulo 3 and  $n$  is  an integer satisfying $n^2\congr D\mod 36$.
The  inclusions $\ca P_{d_i}^{\gamma_2}\subset\ca P_{d_i} $ for $i = 1,2 $  induce a well-defined map
$$
\ca S^{\gamma_2}_{d_1,d_2,n}\longrightarrow\ca S_{d_1,d_2,n} ,\leqno(\eqnumber)
$$
where $S_{d_1, d_2, n} = S_{d_1, d_2, n} ^1 $ is defined as in proposition \eqn2 above.
To prove that this map is a bijection we apply lemma \thmn5 to 
 the  $\Gamma^3$-set $Y = \{(\tau_1, \tau_2)\in\ca P^{\gamma_2}_{d_1}
\times\ca P^{\gamma_2}_{d_2}: B(\tau_1,\tau_2)=n\}$
 and the $\PSL_2(\bo Z)$-set  $X = \{(\tau_1, \tau_2)\in\ca P_{d_1}\times
\ca P_{d_2} :B(\tau_1,\tau_2)=n\}$.
Let $I=\{T^k:k=0,1,2\}$ be a complete set of left coset representatives of 
$\Gamma^3$ in $\PSL_2(\bo Z)$.
As we saw in the proof of proposition \thmn6, the translation $T$ changes the residue class of $b_{\tau_i}$ modulo 3.
Therefore the sets $MY$ with $M\in I$ are disjoint, and for each $(\tau_1,\tau_2)\in X$ there exists
 $M\in I$ such that $b_{M\tau_1}$  is divisible by 3.
Using the assumption on $n$ and the fact that $d_1$ and $d_2$ are congruent to 2 modulo 3
we find by  \eqn2 that $b_{M\tau_2}$ is also divisible by 3.
Hence $M(\tau_1,\tau_2)\in Y$ and  $X$ is the disjoint union of $MY$ with $M\in I$.
It now follows from lemma \thmn5 that \eqn1 is a bijection.
Proposition \thmn1 yields the equality $\rho^{\gamma_2}(n) =\sum_d \varepsilon(d) $ where $d$ 
ranges over the positive divisors of ${1\over 4} (n^2-D)$.
If we sort these divisors by their 3-adic valuation and use the equality 
$\varepsilon(3)=-1$ we obtain
$$
\rho^{\gamma_2}(n)=\cases
\sum_{d|{n^2-D\over 36}} \varepsilon(d) &\hbox{if }n^2\congr D\mod 
36;\cr
0&\hbox{otherwise,}\cr
\endcases\leqno(\eqnumber)
$$
\smallskip \noindent 
Next we sum the Eisenstein series $E_{\P {\Gamma^3}}(z,s)$ over  $z\in\Gamma ^3\back\ca P_{d}^{\gamma_2}$,
where $d$ is equal to~$d_1$ or~$d_2$. 
The stabilizer $\P {\Gamma^3_\infty} $ is generated by  $\ma 1301$ so the scaling transformation at $\infty$ 
is equal to $\sigma_\infty =\ma {\sqrt{3}}00{1\over \sqrt{3}}$.
With $\Gamma =\PSL_2(\bo Z)$ the inclusion $\P{\Gamma^3}\subset \Gamma$ induces a bijection
$\P{\Gamma^3_\infty}\backslash\P{\Gamma^3}\rightarrow \Gamma_\infty\backslash\Gamma$
by lemma \thmn5 and we find the equality $E_{\P {\Gamma^3}} (z,s)=3^{-s}E_\Gamma(z,s)$.
The function that we obtain by substituting $z=\tau$ for some $\tau\in\ca P_d$ in the Eisenstein series $E_\Gamma(z,s)$ is essentially the partial zeta function for the ideal class of $\bo Q(\sqrt{d})$  corresponding to $\tau$. Summing over all $\tau\in\Gamma\back\ca P_d$ yields the zeta function of  $\bo Q(\sqrt{d})$.
More precisely we have the following result [\ZA, proposition~3iii].
\bigskip\noindent\bf Proposition\,\thmnumber\sl
Let $\Delta\congr 0,1\mod 4$ be a negative integer and write $\Delta =df^2$ with~$d$
the discriminant of $\bo Q(\sqrt{\Delta})$.
With $\Gamma=\PSL_2(\bo Z)$ and $\ca P_\Delta $ as in \, $\eqn{35} $ we have
$$
\sum_{\tau\in\Gamma\backslash\ca P_\Delta} E_{\Gamma}(\tau,s)=
{w\over 2}\Big({|\Delta|\over 4}\Big)^{s/2}
{\zeta(s)\over\zeta(2s)} L(s,d)\sum_{m|f}\mu(m){d\legendre m}m^{-s}\sigma_{1-2s}\Big({f\over m}\Big),
$$
where  $w$ is the number of roots of unity in $\bo Q(\sqrt{d})$,  $\zeta(s)$ denotes the Riemann zeta function, $L(s,d)=\sum_{n=1}^\infty{d\legendre n}n^{-s}$ and $\sigma_a(n)=\sum_{m|n} m^a$.
\bigskip\noindent\rm
If we apply this proposition with fundamental discriminant $\Delta =d\congr 2\mod 3 $ and use the bijection
$\P {\Gamma^3}\backslash\ca P_d^{\gamma_2}\rightarrow\Gamma\backslash\ca P_d$, 
which follows from lemma  \thmn6, we obtain
$$
\sum_{\tau\in\Gamma^3\backslash\ca P_d^{\gamma_2}} E_{\P {\Gamma^3}}(\tau,s)=
3^{-s}\hskip -0.2cm\sum_{\tau\in\Gamma\backslash\ca P_d} E_{\Gamma}(\tau,s)=
3^{-s}{w\over 2}\Big({|d|\over 4}\Big)^{s/2}{\zeta(s)\over\zeta(2s)} L(s,d).\leqno(\eqnumber)
$$
Finally we have to calculate the function $\phi_{\P{\Gamma^3}} (s)$.
Fix a positive real number $c$ and let $\sigma_\infty =\ma {\sqrt{3}}00{1\over\sqrt{3}}$ be the scaling transformation at 
 $\infty$.  According to definition  \eqn{10}, we need to compute the number
$\#\left\{d\mod c:\ma **{c\over 3}d \in\P{\Gamma^3}\right\}$. 
For this number to be nonzero $c$ should be an integer multiple of 3.
It follows from the description of $\P{\Gamma^3}$ in section 2 that if ${c\over 3}$ and $d$ are two 
relatively prime integers then  $\P {\Gamma^3} $ contains an element of the form $\ma **{c\over 3}d$.
With $\varphi$ the Euler $\varphi$-function we find
$$
\phi_{\P {\Gamma^3}} (s)=
\sqrt\pi{\Gamma(s-{1\over 2})\over \Gamma(s)}\sum_{c\in 3\bo Z_{>0}} 3 \varphi\big ({c\over 3}\big)c^{-2s} =
3^{1-2s}\sqrt\pi{\Gamma(s-{1\over 2})\over \Gamma(s)}{\zeta(2s-1)\over\zeta(2s)}.\leqno(\eqnumber)
$$
\bigskip\noindent\bf Theorem\,\thmnumber\sl
For $i=1,2$ let $d_i$ and $\alpha_i$ be as in theorem $4$, set $h_i ^\prime ={2\over w_i}h_i $  and
define for $n\in\bo Z_{>0}$ the integer $\rho^{\gamma_2}(n)$ by \, $\eqn3$. 
Then the following formula holds
$$
\log \N\big(\gamma_2(\alpha_1),\gamma_2(\alpha_2)\big)^{8\over w_1w_2}\hskip-0.00cm=
\hskip-0.05cm
\lim_{s\rightarrow 1\atop  s\in\bo H_1}\hskip-0.0cm\left[
{4h_1^\prime h_2^\prime\over s-1} 
-2\hskip -0.15cm\sum_{n>\sqrt{D}}\hskip-0.07cm \rho^{\gamma_2}(n) Q_{s-1}\Big ({n\over\sqrt{D}}\Big)
\right] +4h_1^\prime h_2^\prime C
$$
with
$C={1\over 2}\log D + {L^\prime\over L} (1,d_1)
 + {L^\prime\over L} (1,d_2)-2{\zeta^\prime\over\zeta}(2)-2.
$
\bigskip \noindent\bf Proof.
\rm We need the following Laurent expansions at $s = 1$: 
$$\eqalign{
\zeta(2s)^{-1} =&\textstyle{6\over\pi^2} -{72\zeta^\prime(2)\over\pi^4}(s-1) +O\big((s-1)^2\big) \cr
\textstyle\zeta(s)L(s,d) =&\textstyle{L(1,d)\over s-1} +\gamma L(1,d) +L^\prime(1,d) +O(s-1) \cr
\textstyle{\Gamma(s-\half)\over \Gamma(s)} =&\textstyle\sqrt{\pi} -\sqrt{\pi}\log 4\ (s-1) +O\big ((s-1)^2\big) \cr
\textstyle\zeta(2s-1)=&\textstyle{1\over 2(s-1)} +\gamma+O(s-1) \cr
}$$
where $\gamma$ denotes Euler's constant and $d$ is a negative fundamental quadratic discriminant.
For a derivation of the third expansion above we refer to~[\LA, page 272].
 According to the analytic class number formula we have $L(1,d)={2\over w}{\pi h{|d|}^{-\half}}$ with $h$ and $w$
the class number and the number of roots of unity of the quadratic fields of discriminant $d$. 
The theorem follows from theorem  \thmn4, equations  \eqn{34}, \eqn6, \eqn3, \eqn2, \eqn1 and the expansions above.
\hfill $\qed $
\bigskip\noindent\rm
Next we fix two negative fundamental discriminants $d_1$ and $d_2$ that are relatively prime and congruent
to 1 modulo 8.
In particular both $w_1$ and $w_2$ are equal to~2.
We want to obtain similar results for $\log\N\big(f({-1+\sqrt{d_1}\over 2}),f({-1+\sqrt{d_2}\over 2})\big)$,
where~$f$ denotes one of the functions $\omega$ or $\omega_2$.
According to formula  \eqn{27}  we need to apply theorem \thmn4 for the group  $\overline{\Gamma}_0(2)$ and its normalised principal modulus $-2^{12}/\omega_2$ 
and sum over $z_i\in\overline{\Gamma}_0(2)\back\ca P_{d_i}^f$ ($i =1,2$).

For $f\in\{\omega,\omega_2\}$ the summation of the Green function results in a formula analogous to  \eqn6: 
$$
\sum_{i=1,2,\tau_i\in\overline{\Gamma}_0(2)\backslash\ca P_{d_i}^f}
G_{\overline{\Gamma}_0(2),s}(\tau_1,\tau_2)=
-2\hskip -0.1cm\sum_{n>\sqrt{D}}\rho^f(n)Q_{s-1}(\textstyle{n\over\sqrt{D}})
$$
where $\rho^f(n)$ is the cardinality of the set
$$
\ca S^f_{d_1,d_2,n}=\overline{\Gamma}_0(2)\backslash\{(\tau_1, 
\tau_2)\in\ca P^f_{d_1}\times\ca P^f_{d_2}: B(\tau_1,\tau_2)=n\}.
$$ 
Because of the equalities $\ca P_{d_i} ^{\omega} =\ca P_{d_i} ^2 $ $(i=1,2)$ which we proved in proposition \thmn6, we can apply proposition \thmn3 
with $N =2 $ to find
$$
\rho^\omega(n)=
\cases
2\sum_{d|{n^2-D\over 8}}\varepsilon(d)&\text{\rm  if }n^2\congr D\mod 8;\cr
0&\text{\rm  otherwise.}\cr
\endcases\leqno(\eqnumber)
$$
The computation of $\rho^{\omega_2}(n)$ is more involved.
For $\tau\in\ca P^{\omega_2}_d$ the coefficient $a_\tau $ is odd by the definition in proposition \thmn6.
If the discriminant $d$ is congruent to 1 modulo 8, this implies that~$b_\tau $ is odd and that $c_\tau$
is even.  
Using \eqn6 we find that if $\rho^{\omega_2}(n)$ is nonzero, then~$n$ satisfies $n^2\congr D\mod 16 $.
Therefore we fix a positive integer  $n$ satisfying this congruence.
In the following we suppress the dependence on $d_1,d_2, n $ from the notation.
For positive integers $N, k_1,k_2 $ and $m$ a positive divisor of $N$ we define
$$
\ca A_N(k_1,k_2 ;m) =\overline {\Gamma}_0(N)\back\{(\tau_1,\tau_2)\in\ca P_{d_1}\times\ca P_{d_2}: B(\tau_1,\tau_2) 
= n,a_1\, {\buildrel m\over =}\, k_1,a_2\,{\buildrel m\over =}\, k_2\}
$$
where ${\buildrel m\over =}$ denotes congruence modulo $m$.
This set is well defined because the greatest common divisor of $N$ and $a_\tau$ is 
invariant under the $\overline {\Gamma}_0(N)$-action on $\tau\in\ca P_d$.
To calculate $\rho^{\omega_2}(n) =\#\ca A_2(1,1;2)$ we need the following lemma.
\bigskip\noindent\bf Lemma\,\thmnumber\sl Assume that $d_1$ and $d_2$ are squares modulo $4N$. 
\smallskip
\item{a.} $\#\ca A_N(0,0;N) =2^t\sum_{d|{n^2-D\over 4N}}\varepsilon(d)$, with $t$ the number of prime divisors of  $N$;
\item{b.}  $\#\ca A_N(k_1,k_2;m) ={N\over m}\prod_{p\mid N,p\nmid m}(1+1/p)\#\ca A_m(k_1, k_2;m)$ for $m$ dividing $N$;
\item{c.}  $\#\ca A_4(0,2;4) =\#\ca A_2(0,1;2)$ and $\#\ca A_4(2,0;4)=\#\ca A_2(1,0;2)$;
\item{d.}  $\#\ca A_4(2,2;4) =2\#\ca A_2(1,1;2)$.  
\bigskip \noindent\bf Proof.  
\rm The first formula is the content of proposition \thmn4.
Each of the $\overline {\Gamma}_0(m)$-orbits in $\ca A_m(k_1,k_2;m)$ is the union of $[\overline {\Gamma}_0(m):\overline {\Gamma}_0(N)]$
 $\overline {\Gamma}_0(N)$-orbits.
As $\overline {\Gamma}_0(N)$ is of index $N \prod_{p |N}(1+1/p)$ in $\PSL_2(\bo Z)$, lemma \thmn1b follows.

For $k_1,k_2\in\{0,1\}$ define the set
$$
Y(k_1,k_2) =\{(\tau_1,\tau_2)\in\ca P_{d_1}\times\ca P_{d_2}: B(\tau_1,\tau_2) = n,a_1\,{\buildrel 2\over =}\,k_1, a_2\,{\buildrel 2\over =}\,k_2, c_1\,{\buildrel 2\over =}\, c_2\,{\buildrel 2\over =}\, 0\}.
$$
The greatest common divisors $\gcd(a_\tau,2)$ and $\gcd(c_\tau,2)$ are invariant under the 
 $\overline{\Gamma}(2)$-action on $\tau\in\ca P_d$.
Hence $\overline{\Gamma}(2)$ acts diagonally on $Y(k_1,k_2)$.
For $i =1,2$ the maps $\tau_i \mapsto 2\tau_i$, or on coefficients  $[a_i, b_i,c_i] \mapsto [{a_i\over 2},b_i,2c_i]$,
induce a map
$$\ca A_4(2k_1,2k_2;4)\rightarrow \overline{\Gamma}(2)\back Y(k_1,k_2). $$
It follows from the equality
$\ma 2001 \overline {\Gamma}_0(4) =\overline{\Gamma}(2)\ma 2001$ that this map is a bijection.
If $\{k_1,k_2\} =\{0,1\} $ we can apply lemma \thmn8 to find that $\overline{\Gamma}(2)\back Y(k_1,k_2) $  
has the same cardinality as $\ca A_2(k_1,k_2;2)$ and lemma \thmn1c follows.
As the discriminants are congruent to 1 modulo 8 by assumption, the conditions on $c_1,c_2$
in the definition of $Y(1,1)$ follow from those on $a_1,a_2$.
In particular we have the equality $\ca A_2(1,1;2) =\overline{\Gamma}_0(2)\back Y(1,1)$. 
Each of the $\overline{\Gamma}_0(2)$-orbits in this set is the union of two $\overline{\Gamma}(2)$-orbits
and the last formula of the lemma follows.\hfill$\qed$
\bigskip\noindent\rm 
By counting the complement of $\ca A_2(1,1; 2) $ inside $\ca A_2(0,0;1)$ we find the equality
$$
\#\ca A_2(1,1;2) =\#\ca A_2(0,0;1)-\#\ca A_2(0,0;2) - \#\ca A_2(0,1;2) -\#\ca A_2(1,0;2).\leqno(\eqnumber)
$$
To compute the first two terms on the right hand side we apply lemma \thmn1a,b and obtain
$\ca A_2(0,0;1) =3\ca A_1(0,0;1) =3\sum_{d|{\ba n\over 4}} \varepsilon(d)$ and $\ca A_2(0,0;2) =
2\sum_{d|{\ba n\over 8}} \varepsilon(d)$ with $\ba n=n^2-D $.
For the last two terms we first use lemma \thmn1c and then count the complement of 
$\ca A_4(0,2;4)\cup \ca A_4(2,0;4)$ inside $\ca A_4(0,0;2)$: 
$$
\#\ca A_2(0,1;2) +\#\ca A_2(1,0;2) =\#\ca A_4(0,0;2) -\#\ca A_4(0,0;4) -\#\ca A_4(2,2;4).
$$
As above we apply lemma \thmn1a,b to calculate the first two terms on the right of this equation.
By lemma \thmn1d the last term is equal to $2\#\ca A_4(1,1;2)$, twice the number we are trying to calculate.
Substituting all these terms into equation \eqn1 yields
$$
\#\ca A_2(1,1; 2) = -3\sum_{d |{\ba n\over 4}}\varepsilon(d) +2\sum_{d | {\ba n\over 8}}\varepsilon(d)
+4\sum_{d | {\ba n\over 8}}\varepsilon(d) -2\sum_{d | {\ba n\over 16}}\varepsilon(d).
$$
Using that $\varepsilon$ is multiplicative and $\varepsilon(2) = 1 $ we conclude
$$
\rho^ {\omega_2} (n)=
\cases
\sum_{d|{n^2-D\over 16}}\varepsilon(d)&\text{\rm  if }n^2\congr D\mod 16;\cr
0&\text{\rm  otherwise.}\cr
\endcases\leqno(\eqnumber)
$$
Let $d$ denote one of the discriminants $d_1$ for $d_2$  and set $\Gamma=\PSL_2(\bo Z)$. 
To calculate the summation of $E_{\overline{\Gamma}_0(2)}(z,s)$ over $z\in\overline{\Gamma}_0(2)\back\ca P ^\omega_d $ we use proposition \thmn3 and the relation
$$
E_{\overline{\Gamma}_0(2)}(z,s)=
{2^s\over2^{2s}-1}\big( E_\Gamma(2z,s)-2^{-s}E_\Gamma(z,s)\big)
$$
between the Eisenstein series for $\Gamma$ and $\overline{\Gamma}_0(2)$ [\GZt, \S 2 (2.16)].
In the proof of proposition \thmn7 we showed that the map $\overline{\Gamma}_0(2)\back\ca P ^\omega_d\rightarrow\Gamma\back\ca P_d $ induced by the inclusion $\ca P ^\omega_d\subset\ca P_d $ is 2 to 1.
In a similar way one proves that the map $\overline{\Gamma}_0(2)\back\ca P ^\omega_d\rightarrow\Gamma\back\ca P_d $ sending $\tau$ to $2\tau $ is also 2 to 1
and we obtain
$$\eqalign {
\sum_{\tau\in\overline{\Gamma}_0(2)\backslash\ca P_d^{\omega}}E_{\overline{\Gamma}_0(2)}(\tau,s) & =
{2^{s}\over2^{2s}-1}\sum_{\tau\in\Gamma\backslash\ca P_d}
\big(2E_\Gamma(\tau,s)-2^{1-s}E_\Gamma(\tau,s)\big) \cr
& ={2\over 2^s+1}\Big({|d|\over 4}\Big)^{s/2}{\zeta(s)\over\zeta(2s)} L(s,d)\cr
} $$ 
To obtain a similar result for $\omega_2$ we let 
$Y =\{\tau\in X:2\nmid a_\tau,2 ||b_\tau\} $ be the image of $\ca P^{\omega_2}_d$ in
$X =\left\{\tau\in\ca P_{4d}:\gcd(a_\tau,b_\tau,c_\tau) =1\right\}$
under the map $\tau\mapsto2\tau $.
If we define 
$\overline{\Gamma} ^0(2)=\ma 2001\overline{\Gamma}_0(2) {\ma 2001}^{-1}$ we conclude that the first map below is a bijection.
$$
\overline {\Gamma}_0 (2)\back\ca P ^{\omega_2}_d\ 
{\buildrel \tau\mapsto 2\tau \over \longrightarrow}\ 
\overline{\Gamma} ^0(2)\back Y\ 
\longrightarrow\ 
\Gamma\back X
$$  
By lemma \thmn8 the second map, induced by the inclusion $Y\subset X$, is also a bijection.
Using that the complement of $X$ in $\ca P_{4d}$ is equal to $\ca P_d$ 
and the bijection \eqn{32} between
$\overline{\Gamma}_0(2)\back\ca P ^ {\omega_2}_d $ and $\Gamma\back\ca P_d $
we find
$$\eqalign{
\sum_{\tau\in\overline{\Gamma}_0(2)\backslash\ca P_d^ {\omega_2}}\hskip -0.4 cm E_{\overline{\Gamma}_0(2)}(\tau,s) & =
{2^{s}\over2^{2s}-1}\Big[\hskip -0.1cm
\sum_{\tau\in\Gamma\backslash\ca P_{4d}}\hskip -0.3 cm E_\Gamma(\tau,s)-\hskip -0.2 cm
\sum_{\tau\in\Gamma\backslash\ca P_{d}}\hskip -0.25 cm E_\Gamma(\tau,s)-
2^{-s}\hskip -0.25 cm\sum\limits_{\tau\in\Gamma\backslash\ca P_{d}}\hskip -0.25 cm E_\Gamma(\tau,s)\Big]\cr
& ={2^s-1\over 2^s+1}\Big({|d|\over 4}\Big)^{s/2}{\zeta(s)\over\zeta(2s)} L(s,d).\cr
} $$
The scaling transformation $\sigma_\infty$ for the group $\overline{\Gamma}_0(2)$ is the identity,
hence definition \eqn{14} yields 
$$
\phi_{\overline{\Gamma}_0(2)}(s)=
\sqrt\pi{\Gamma(s-{1\over 2})\over \Gamma(s)}\sum_{s\in2 \bo Z_{>0}}{\varphi(c)\over c ^{2s}} =
 {\sqrt{\pi}\over 4^s-1}{\Gamma(s-{1\over 2})\over \Gamma(s)}{\zeta(2s-1)\over\zeta(2s)}.
$$
Combining the above formulas as in the proof of theorem \thmn2 we find the following result.
\bigskip\noindent\bf Theorem\,\thmnumber\sl
For $i=1,2$ let $d_i$ and $\alpha_i$  be as in theorem $5$.
Define for $n\in\bo Z_{>0}$ the integers $\rho^\omega(n)$ and $\rho^{\omega_2}(n)$ 
by  $\eqn3$ and  $\eqn1$ respectively.
Then the following formulas hold:
$$\displaylines{
\textstyle
\log \N\big(\omega(\alpha_2), \omega(\alpha_2)\big)\hskip -0.0 cm =
\hskip -0.0 cm\lim\limits_{\scriptstyle s\rightarrow 1\atop \scriptstyle s\in\bo H_1}
\hskip -0.1 cm\Big [
{8h_1 h_2\over s-1}
-\hskip-0.2cm\sum\limits_{n>\sqrt{D}}
\textstyle\hskip -0.1 cm\rho^{\omega}(n)Q_{s-1}\big({n\over\sqrt{D}}\big)
\Big] 
+8h_1h_2(C+\textstyle{4\log2\over 3})\cr
\textstyle
\log \N\big(\omega_2(\alpha_2), \omega_2(\alpha_2)\big)\hskip -0.0 cm =
\hskip -0.00 cm\lim\limits_{\scriptstyle s\rightarrow 1\atop \scriptstyle s\in\bo H_1}
\hskip -0.1cm\Big [
{2h_1 h_2\over s-1}
-\hskip-0.2cm\sum\limits_{n>\sqrt{D}}
\textstyle\hskip -0.1cm\rho^{\omega_2}(n)Q_{s-1}\big({n\over\sqrt{D}}\big)
\Big] 
+2h_1h_2(C-\textstyle{2\log2\over3})\cr
}$$
with $C$ as in theorem $13$.

\head 4. A family of non-holomorphic Hilbert modular forms and holomorphic projection
\endhead
\noindent\rm
In order to finish the proof of the theorems in the introduction we need to compute the limit
$$
\lim_{s\rightarrow 1\atop s\in\bo H_1}\left[
{4h_1^\prime h_2^\prime\over s-1}
-2\hskip -0.1cm\sum_{n>\sqrt{D}}\sum_{d |{n ^2-D\over 36}}
\hskip -0.1 cm\varepsilon(d)
Q_{s-1}\Big ({n\over\sqrt{D}}\Big) \right] 
\leqno(\eqnumber)
$$
which occurs in theorem   \thmn3 and the limits occurring in theorem  \thmn1.
Gross, Kohnen and Zagier proved [\GZd, \S 3] that limits similar to the one above are `almost' equal to the first Fourier coefficient
of certain cusp forms.
As their proof does not cover the above limit we extend their result to a slightly larger class of limits.
In this section, which is independent of the previous ones, we follow Gross and Zagier and study a family of non-holomorphic Hilbert modular forms.
Via holomorphic projection we obtain a family of holomorphic cusp forms of weight 2 on congruence subgroups of $\SL_2(\bo Z)$.
The main result, and the only one which we will use in the sequel, is theorem  \thmn{-4} below which states that the first Fourier coefficient of each of these cusp forms 
is up to a `simple expression' equal to a limit like \eqn1.
In the next section we concentrate on those cusp forms in the family for which the involved congruence subgroup has genus zero.
The corresponding cusp form is then identically zero and hence \eqn1 is equal to the `simple expression', 
which completes the proof of the main theorems.

To indicate the relation between our main theorems and the Fourier coefficients of modular forms we recall
the motivation given by Gross and Zagier on page 214 of [\GZe].
Let  $K$  be a real quadratic field and denote its Dedekind zeta function by~$\zeta_K$.
For $k$ a positive even integer Siegel expressed $\zeta_K(1-k)$ as a finite sum of norms of certain ideals [\SI].
For $k = 2,4$ these expressions are given by the first equality in the following formula:
$$
30k\zeta_K(-k+1) =\sum_{{\nu\in\go d^{-1}\atop\nu \gg 0}\atop \Tr(\nu) =1}\sum_{\go n | (\nu)\go d}N_{K/\bo Q}(\go n)^{k-1} =\sum_{x ^2 <D\atop x ^2\congr D\mod 4}\sum_{n |{D- x ^2\over 4}}n^{k-1}\leqno(\eqnumber)
$$
where $\go d =(\sqrt{D})$ is the different of  $K$ and where we use $\nu \gg 0$ to denote 
that $\nu$ is totally positive.
The second equality above follows by noting that the totally positive $\nu\in\go d^{-1}$ of trace 1 are of the form
$\nu={x+ \sqrt{D}\over 2\sqrt{D}}$ for some integer $x$ which satisfies  $x\congr D\mod 2$ and $x^2<D$.
Siegel's proof of the first equality uses the Hecke-Eisenstein series $E_{K,k}(z_1,z_2)$,
a holomorphic Hilbert modular form of weight $k$ on $\SL_2(\ca O)$ with  $\ca O$ the ring of integers of $K$.
To define this series we fix and embedding $K\rightarrow \bo C$, denote the conjugate of $x\in K$
by $\ba x$ and let $C$ be the ideal class group of $K$.   
For  $k\in2\bo Z_{>0} $ and $(z_1,z_2)\in\bo H\times\bo H$ we define
$$
E_{K,k}(z_1, z_2)=\hskip -0.1cm\sum_{[\go a]\in C} 
\sum_{\scriptstyle (m,n)\in (\go a\times\go a)/\ca O^*\atop\scriptstyle (m,n)\not=(0,0)}
\hskip-0.1cm
{N_{K/\bo Q}(\go a)^{k}\over (m z_1+n)^k (\ba m z_2+\ba n)^k},
$$
where in the case $k =2$ the sum is computed by Hecke summation [\SI, page~93].
The restriction  $F_{K, k} (z) =E_{K, k} (z,z)$ to the diagonal is an ordinary modular form of weight $2k$ on $\SL_2(\bo Z)$.
Its zeroth and first Fourier coefficient are equal to~$\zeta_K(k)$ and
$({2\pi\over (k-1)!})^{2k}D^{{\half}-k}$ times the middle term of  \eqn1, respectively.
In case~$k$ is equal to 2 or 4 the vector space of modular forms of weight $2k$ on $\SL_2(\bo Z)$ is 
one-dimensional and hence $F_{K,k}(z)$ is a multiple of the ordinary Eisenstein series of weight $2k$.
Using the functional equation of~$\zeta_K$ we obtain equation \eqn1.

The right hand side of \eqn1 resembles the right hand side of the formula for $J(d_1,d_2)$ (see theorem 1):
$$
-\log |J(d_1,d_2) |^{8\over w_1w_2} =\sum_{x^2 < D\atop x^2\congr D\mod 4}\sum_{n |{D- x^2\over 4}}\varepsilon(n)\log n.\leqno(\eqnumber)
$$
In order to prove this formula Gross and Zagier adapt the definition of $F_{K,k}(z)$
guided by the distinction
between \eqn1 and \eqn2 and study its Fourier coefficients.
As our theorems are similar to \eqn1 we will do the same.

For the rest of this section we assume that the discriminant $D$ of $K$ can be written in the
form $D =d_1d_2 $ with $d_1$ and $d_2$  two fixed relatively prime negative fundamental discriminants.
To account for $\varepsilon$ in \eqn1 we introduce the genus character $\chi:C^+\rightarrow \{\pm 1\}$ of the strict ideal class group
 $C^+$ of $K$ corresponding to the decomposition $D=d_1\cdot d_2$.
If we identify $C^+ $ via the Artin map with the Galois group of the strict Hilbert class field
 $H^+ $ over $K$,  the character $\chi$ has kernel $\Gal(H^+/\bo Q(\sqrt{d_1},\sqrt{d_2}))$.
We extend $\chi$ to the idealgroup of $K$ and find $\chi(\go p)=1$ for inert primes $\go p$ 
and  $\chi(\go p)=\varepsilon(N_{K/\bo Q}(\go p))$ for the other primes $\go p$.
As $\bo Q(\sqrt{d_1},\sqrt{d_2})/K$ is ramified at infinity, the character $\chi$ does not factor through
the ideal class group $C$  and hence $\chi\big ((\lambda )\big)=\sign\big (N_{K/\bo Q} (\lambda)\big) $ for $\lambda\in K^*$.

To transform the term $n^{k-1}$ of \eqn2  into the term $\log n$ of  \eqn1 it seems natural
`to differentiate $F_{K,k}(z)$  with respect to  $k$ and substitute $k = 1$'. 
However the holomorphic Hecke-Eisenstein series is only defined for positive even integers $k$.
To resolve this difficulty we introduce for complex $s$  the factor ${{y_1y_2}^s\over |mz_1+n |^{2s} | \ba mz_1+\ba nz_2 |^{2s}}$ into the non-convergent series  $E_{K,1}(z_1,z_2)$ (`Hecke's trick'),
restrict to the diagonal  $z_1 =z_2 $  and take the derivative with respect to $s$ at $s = 0$. 

Fix a primitive ideal $\go n\subset\ca O$  of norm $N$ such that $\chi(\go n)=\varepsilon(N) =1$.
(We do not assume $\chi$ to be trivial on all divisors of $\go n$ as in [\GZd].)
For complex $s$ with $\Re(s)>\half$ we define the following non-holomorphic Hecke-Eisenstein
series on $\bo H\times\bo H$: 
$$
E_{\go n,s}(z_1, z_2)=\hskip -0.1cm\sum_{[\go a]\in C} 
\sum_{\scriptstyle (m,n)\in (\go a\go n\times\go a)/\ca O^*\atop\scriptstyle (m,n)\not=(0,0)}
\hskip-0.1cm
{\chi(\go a)N_{K/\bo Q}(\go a)^{1+2s} (y_1y_2)^s\over (m z_1+n)(\ba m z_2+\ba n)| m 
z_1+n|^{2s}| \ba m z_2+\ba n|^{2s}},
$$
with $y_1 =\Im(z_1)$ and $y_2 =\Im(z_2)$.  
Although $\chi$ is a character of the strict ideal class group the summation over $[\go a]\in C$ is well defined;
if we replace $\go a$ by $\lambda\go a$ for some $\lambda\in K^*$ the summand does not change
because of the equality $\chi\big ((\lambda )\big)=\sign\big (N_{K/\bo Q} (\lambda)\big) $.
As $D$ is the product of two negative fundamental discriminants
the elements of the unit group $\ca O^*$ have norm 1
and the inner summation over the $\ca O^*$-orbits is also well defined.
For $(z_1,z_2) $ in some fixed compact subset of $\bo H\times\bo H$ the estimates 
$|mz_i+n | =O(\sqrt{m^2+n^2}) $ for $i =1,2$ lead to $E_{\go n,s}( z_1, z_2) =O\big(\zeta_{K(\sqrt{-1})}(s+\half)\big)$.
The above sum is therefore absolutely and locally uniformly convergent for 
$( z_1, z_2,s)\in\bo H\times\bo H\times\{s\in\bo C:\Re(s) >\half\}$.
In particular $E_{\go n,s}( z_1, z_2)$ is analytic as function of $s$ for $\Re(s) >\half$.
In fact it has a meromorphic continuation to all of $\bo C$ as we will now show.

Fix $s\in\bo C $ with $\Re(s) >\half$. 
Using the absolute convergence one easily finds that  $E_{\go n,s}( z_1, z_2) $ transforms like a Hilbert modular form of weight~1:
$$
E_{\go n,s}\Big( {az_1+b\over cz_1+d},{\ba a z_2+\ba b\over \ba c z_2+\ba d} \Big)=
(cz_1+d)(\ba c z_1 +\ba d)E_{\go n,s}(z_1,z_2)\leqno(\eqnumber)
$$
for all $\ma abcd\in\PSL(\ca O,\go n)$ 
with
$$
\PSL(\ca O,\go n) =\left\{\ma abcd\in\PSL_2(K): a,d\in\ca O, c\in\go n \text{\rm \ and\ }b\in\go n^{-1}\right\}.
$$
In particular we find $E_{\go n,s}(z_1+\lambda,z_2+{\ba \lambda}) =E_{\go n,s}(z_1,z_2) $ for each $\lambda\in\go n^{-1}$ and the Hilbert modular form
can be viewed as a function of $\big(y_1, y_2,(x_1,x_2)\big)\in\bo R_{>0}\times\bo R_{>0}\times\bo R ^2/\go n^{-1}$. 
Here and below we identify any fractional $K$-ideal $\go a$ with the lattice $\{(\lambda,\ba\lambda):\lambda\in\go a\}\subset\bo R^2$. 
The compact group $\bo R ^2/\go n^{-1}$ has Lebesque measure $\sqrt{D}N^{-1}$ and its characters
are of the form $(x_1,x_2)\mapsto e^{-2\pi i(\nu x_1+\ba \nu x_2)}$ with $\nu$ in the dual lattice
$$
\left\{\nu\in K ^*:\Tr_{K/\bo Q}(\nu \go n^{-1})\subset\bo Z\right\} =\go n\go d^{-1}.
$$
Hence $E_{\go n,s}(z_1,z_2) $ has a Fourier expansion of the form
$$
E_{\go n,s}(z_1,z_2)=\sum_{\nu\in\go n\go d^{-1}} c_{\nu,s}(y_1,y_2)
e^{2\pi i (\nu x_1+\ba \nu x_2)}\leqno(\eqnumber)
$$
with
$$
c_{\nu,s}(y_1,y_2) ={N\over\sqrt{D}}\int_{\bo R ^2/\go n^{-1}}E_{\go n,s}(z_1,z_2)e^{-2\pi i(\nu x_1+\ba \nu x_2)}dx_1dx_2\leqno(\eqnumber)
$$
where $x_i=\Re(z_i)$ and $y_i=\Im(z_i)$ for $i=1,2$ and $dx_1dx_2$ denotes the standard Lebesque measure.
In the proposition below we will calculate these Fourier coefficients.  
As it will turn out they can be meromorphically continued to the whole 
 $s$-plane which gives us the continuation of $E_{\go n,s}(z_1,z_2)$ via \eqn2. 

We need to introduce the following functions.
For $\go a$ a nonzero integral ideal of~$\ca O$ and $s\in\bo C $ define the function 
$$
\sigma_{s,\chi}(\go a)=\sum_{\go b|\go a}\chi(\go b)N_{K/\bo Q}(\go b)^s\leqno(\eqnumber)
$$
which is clearly analytic in $s$.
This definition also makes sense in case $\go a$ is the zero ideal if we let $\go b$ range over the nonzero integral ideals of $\ca O$ and restrict $s$ to the half plane $\Re(s) < -1 $.
For $\Re(s) >1 $ we then have the equality $\sigma_{-s,\chi}\big((0)\big) = L_K(s,\chi) $ 
with
$$
L_K(s,\chi)=\sum_{0\not=\go a\subset\ca O}\chi(\go a)N_{K/\bo Q}(\go a)^{-s}
$$
the Hecke $L$-function associated to the character $\chi$.  
Finally for $s\in\bo C,\Re(s) >0 $ and $t\in\bo R $ we define the function 
$$
\Phi_s(t)=\int_{-\infty}^\infty{e^{-2\pi ixt}\over(x+i)(x^2+1)^s}dx.\leqno(\eqnumber)
$$
\bigskip\noindent\bf Proposition\,\thmnumber\sl 
For $\Re(s) >\half $  and $\nu\in\go n\go d^{-1}$  the Fourier coefficient  $\eqn3 $ is given~by
$$
{c_{\nu,s}(y_1,y_2) =\cases
L_K(1+2s,\chi)(y_1y_2)^s + {\textstyle\Phi_s(0)^2L_K(2s,\chi)\over\textstyle 
\sqrt{D} N^{2s}(y_1y_2)^{s} }&\hbox{if }\nu=0;\cr
{\textstyle\Phi_s(\nu y_1)\Phi_s(\ba\nu y_2)\sigma_{-2s,\chi}\big((\nu)\go d\go n^{-1}\big)\over\textstyle
\sqrt{D} N^{2s}(y_1y_2)^{s} } &\hbox{if }\nu\not=0.\cr
\endcases}\leqno(\eqnumber) 
$$
\bf Proof. \rm The following computation is a variation of the methods of Hecke [\HEc, \S 3].
See [\GZe, page 214] for a discussion on an error of sign made by Hecke.

We split the summation $(m, n)\in (\go a\go n\times\go a)/\ca O ^* $ in the definition of $E_{\go n,s}(z_1,z_2)$ into
two: one sum over  $m =0 $ and $n\in (\go a -\{0\})/\ca O ^*$ and one sum over $m\in (\go a\go n-\{0\})/\ca O ^*$ and $n\in\go a$.
The contribution of the first summation to the integral in~\eqn4  is
$$
\sum_{[\go a]\in C}\sumprime_{n\in\go a/\ca O ^*} {\chi(\go a)N(\go a)^{1+2s}(y_1y_2) ^s\over N(n) |N(n) |^{2s}}\int_{\bo R ^2/\go n^{-1}}e^{-2\pi i(\nu x_1+\ba \nu x_2)}dx_1dx_2.
$$
Here and below a prime above a summation sign indicates that the summation variable ranges over nonzero  elements.
Furthermore $N(\go a)$ and $N(n)$ are abbreviations of the norms $N_{K/\bo Q} (\go a) $ and $N_{K/\bo Q}(n)$.
The integral above is zero unless $\nu=0$ in which case it equals $\sqrt{D} N^{-1}$.
The map $n\mapsto n\go a^{-1} $ is a bijection between $(\go a -\{0\})/\ca O ^* $ and the set of integral ideals in the ideal class of $\go a^{-1}$.
Using this bijection we conclude that the first summation contributes nothing to $c_{\nu,s} $  for $\nu\ne0$ 
and 
$$
L_K(1+2s,\chi)(y_1y_2) ^s
$$ 
to $c_{0,s}$.

Fix a complete set of representatives  $S$ for $(\go a\go n-\{0\})/\ca O ^* $.
The summation over $m\in (\go a\go n-\{0\})/\ca O ^*$ and $n\in\go a$ yields the following contribution to the integral in~\eqn4:
$$
\sum_{[\go a]\in C}\sum\limits_{m\in S} {\chi(\go a)N(\go a)^{1+2s} \over N(m) |N(m) |^{2s}} 
\hskip-0.05cm\int_{\bo R ^2/\go n ^{-1}}\hskip -0.05cm\sum_{n\in\go a}{(y_1y_2)^se^{-2\pi i(\nu x_1+\ba \nu x_2)}dx_1dx_2\over (z_1+{n\over m})(z_2+{\ba n\over \ba m}) |z_1+{n\over m} |^{2s} |z_2+{\ba n\over \ba m} |^{2s}}.
$$
Fix an element $m\in S$.
For $n\in\go a$ we write  $n=n_1+m n_2$ with $n_1\in\go a/m\go n^{-1}$ and $n_2\in\go n^{-1}$.
Summing over $n_2$ and writing $n$ for $n_1$, the above integral is equal to  
$$
\int_{\bo R ^2}\sum_{n\in\go a/m\go n^{-1}}{(y_1y_2)^se^{-2\pi i(\nu x_1+\ba \nu x_2)}dx_1dx_2\over (z_1+{n\over m})(z_2+{\ba n\over \ba m}) |z_1+{n\over m} |^{2s} |z_2+{\ba n\over \ba m} |^{2s}}
$$
where we used the equality $e^{-2\pi i(\nu n_2+\ba {\nu n_2})} =1$ following from the assumption $\nu\in\go n\go d^{-1}$. 
The change of variables $x_1=y_1\tilde{x}_1-{ n\over m } $ and $x_2 =y_2\tilde{x}_2-{\ba n\over\ba m} $ yields
$$(y_1y_2)^{-s}\Phi_s(\nu y_1)\Phi_s(\ba \nu y_2)
\sum_{n\in\go a/m\go n^{-1}}e^{2\pi i({\nu\over m} n  + {\ba \nu\over\ba m}\ba n)}.
$$
The finite sum above is zero unless ${\nu\over m}\in (\go a\go d)^{-1} $ in which case it equals ${| N(m) |\over 
N(\go a\go n)}$.
Hence the summation over $m\in (\go a\go n-\{0\})/\ca O ^*$ and $n\in\go a$ yields the following contribution to $c_{\nu,s}(y_1,y_2)$: 
$$
{\Phi_s(\nu y_1)\Phi_s(\ba \nu y_2)\over \sqrt{D}(y_1y_2)^{s}}
\sum_{[\go a]\in C}\sumprime_{\scriptstyle m\in\go a\go n/\ca O ^*:\atop\scriptstyle \nu\in m (\go a\go d)^{-1}}
{ \chi(\go a)N(\go a)^{2s} \over N(m)|N(m) |^{2s-1}}.
$$
The elements $m$ in the inner summation are in bijection with the integral ideals in the ideal class of $(\go a\go n)^{-1}$ that contain $(\nu)\go d\go n^{-1}$.
By changing the outer summation $[\go a]\in C$ into $[(\go a\go n)^{-1}]\in C$ we find that the above double sum is equal to $N^{-2s}\sigma_{-2s,\chi} \big((\nu)\go d\go n^{-1}\big)$.
Using the equality $\sigma_{-2s,\chi}((0)) = L_K(2s,\chi) $ this concludes the proof of the proposition.\hfill$\qed$
\bigskip\noindent\bf Corollary\,\thmnumber\sl 
The function $E_{\go n,s}(z_1,z_2) $ has a meromorphic continuation to the entire $s$-plane which is regular at
 $s =0 $. 
Moreover we have $E_{\go n,0}(z_1,z_2) =0$. 

\medskip \noindent\bf Proof. \rm First we recall some properties of the function $\Phi_s(t)$ defined by  \eqn2,
proofs can be found in [\GZt, IV \S 3]. 
For $t =0 $ and $\Re(s) >0 $  we have the equality
$$
\Phi_s(0) = -\sqrt{\pi} i{\Gamma(s+\half)\over\Gamma(s+1)}.\leqno(\eqnumber)
$$
Hence  $\Phi_s(0)$ admits a meromorphic continuation to the whole $s$-plane which is regular at $s =0 $. 
By changing the path of integration one can show that for all nonzero $t\in\bo R $ the function
$\Phi_s(t)$ has an analytic continuation to $s\in\bo C$.
For each compact $V\subset\bo C $ there exist positive constants $c_1,c_2 $ such that
$$
|\Phi_s (t) |\le {c_1\over |t |^{c_2}}e^{-\pi |t|}\text{\rm \ \ for all\ } s\in V,\ t\ne 0, \leqno(\eqnumber)  
$$
and at $s =0$ the function is equal to 
$$
\Phi_0(t) =\cases
-2\pi ie^{-2\pi t} &\text{\rm  if\ }t >0;\cr
-\pi i & \text{\rm if\ }t =0;\cr
0 & \text{\rm  if\ }t <0.
\endcases\leqno(\eqnumber)
$$
Recall that for $\go a$ a nonzero integral ideal the function $\sigma_{s,\chi}(\go a)$ defined by  \eqn6 is analytic on $\bo C$. 
As is well known [\NE, VII \S 8] the Hecke $L$-series  $L_K(s,\chi)$ has an analytic continuation to $s\in\bo C$
and satisfies the functional equation 
$$
D^{s/2}\pi^{-(s+1)}\Gamma(\textstyle{s+1\over 2}) ^2 L_{K}(s,\chi) =D ^{(1-s)/2}\pi^{s-2}\Gamma({2-s\over 2}) ^2 L_K(1-s,\chi).\leqno(\eqnumber)
$$
Hence all the functions occuring on the right hand side of  \eqn5 have  meromorphic continuations to the whole $s$-plane which are regular at $s =0$, and we conclude that the same holds for $c_{\nu,s}(y_1,y_2)$. 
Moreover for nonzero $\nu$ this continuation is analytic in $s$.

Because of \eqn3 the Fourier series  \eqn9 of $E_{\go n,s}(z_1,z_2)$ converges absolutely and uniformly
on compact subsets of $\bo C$ that do not contain the poles of $\Phi_s(0)$.
This proves the continuation as claimed in the corollary. 

To prove that $E_{\go n,s}(z_1,z_2) $ vanishes at  $s =0 $ we show that its Fourier coefficients  \eqn5 are all zero. 
The equality $c_{0,0}(y_1,y_2) =0 $ follows from \eqn4 and the functional equation  \eqn1 for $L_K(s,\chi)$.
If either $\nu$ or $\ba\nu $ is negative,  $c_{\nu,s}(y_1,y_2)$ is zero because of \eqn2.
Now assume $\nu$ is totally positive. As by assumption $\chi$ is trivial on $\go n$  we find
 $\chi((\nu)\go d\go n^{-1}) =\chi(({\nu\over\sqrt{D}})) = -\sign(N(\nu)) <0$.
Hence the contribution of $\go b |(\nu)\go d\go n^{-1} $ to $\sigma_{0,\chi}((\nu)\go d\go n^{-1})$ is cancelled by
the contribution of the complementary divisor  $(\nu)\go d (\go n\go b) ^{-1}$.
We find $\sigma_{0,\chi}((\nu)\go d\go n^{-1})=0$ so that $c_{\nu,0}(y_1,y_2)$ vanishes for totally positive~$\nu$.\hfill $\qed$
\bigskip\noindent As we argued in the beginning of this section we are interested in the derivative of the Hecke-Eisenstein series at $s =0 $ restricted to the diagonal:
$$F_\go n(z)={\sqrt{D}\over 8\pi^2}{\partial\over\partial 
s}E_{\go n,s}(z,z)\big|_{s=0}.$$
According to  \eqn{10} this is a non-holomorphic modular form of weight 2 on $\Gamma_0(N)$.
The Fourier expansion of $E_{\go n,s}(z_1,z_2)$ is locally uniformly convergent in $s$  because of  \eqn3. 
By differentiating  \eqn9  termwise and using  \eqn5, \eqn4, \eqn2 and \eqn1  we obtain the 
Fourier expansion of $F_\go n(z) $: 
$$
F_\go n(z)=\sum_{n = -\infty}^{\infty}\Big(\hskip -0.25cm \sum_{\scriptstyle \nu\in\go n\delta^{-1}\atop\scriptstyle  \Tr(\nu) = n} \hskip -0.1cm c_\nu(y)\,\Big)
\ e^{2\pi inz} \leqno(\eqnumber)
$$
where 
$$
c_\nu(y)=\cases
{\sqrt{D}\over 2\pi^2}(L_K(1,\chi)\log y +\kappa)&\hbox{if }\nu=0;\cr
{\sigma}^\prime_\chi\big( (\nu)\go d\go n^{-1}\big)&\hbox{if }\nu\gg 0;\cr
-{1\over 2}\sigma_{0,\chi}\big( (\nu)\go d\go n^{-1}\big)\Phi(|{\bar\nu}|y)
&\hbox{if }\nu>0>\ba \nu;\cr
-{1\over 2}\sigma_{0,\chi}\big( (\nu)\go d\go n^{-1}\big)\Phi(|\nu|y)
&\hbox{if }\ba \nu>0>\nu;\cr
0&\hbox{if }\nu<<0,\cr
\endcases\leqno(\eqnumber)
$$
with
$$
\displaylines {
\textstyle \kappa =L_K^\prime(1,\chi)+\big (\half\log(DN)-\log \pi-\gamma 
\big)L_K(1,\chi),\cr
\textstyle {\sigma}^\prime_\chi(\go a)= {\partial\over\partial s}\sigma_{s,\chi}(\go a)
\big|_{s=0}=\sum_{\go m|\go a}\chi(\go m)\log N_{K/\bo Q}(\go m)\cr
\text{\rm and} \hfill \cr
\textstyle\Phi(t)={i\over 2\pi}e^{-2\pi t}{\partial\over\partial s}\Phi_s(-t)\big|_{s=0}\hskip 0.5 cm\text{\rm  for\ }t\in\bo R_{>0} .\cr
}
$$
By differentiating under the integral sign one proves the equality 
$$
\Phi(t) =\int_1 ^\infty e^{-4\pi tu}{du\over u}\hskip 0.5 cm\text{\rm  for\ }t >0, \leqno(\eqnumber)
$$
[\GZt, IV \S 3] from which we obtain
$$
\Phi(t) =O(t^{-1}e^{-4\pi t })\hskip 0.5 cm\text{\rm  for\ }t >0,\leqno(\eqnumber)
$$
which guarantees the convergence of  \eqn4.

To see why these Fourier coefficients can be used to compute limits like  \eqn{17}, we have to introduce the notion of holomorphic projection.
Let $S_2(N)=S_2(\Gamma_0 (N))$ be the finite dimensional complex vector space of 
holomorphic cusp forms of weight 2 on $\Gamma_0(N)$.
This vector space is equipped with a Hermitian inner product which for $f,g\in S_2(N)$
is defined by
$$
\langle f, g\rangle =\int_{\Gamma_0(N)\back\bo H}f(z)\overline{g(z)}y^2d\mu,\leqno(\eqnumber)
$$
where the integral is taken over a fundamental domain for the action of $\Gamma_0(N)$ on $\bo H$
and where $d\mu $ denotes the  $\PSL_2(\bo R)$-invariant measure $dxdy\over y ^2$ on $\bo H$  with $z =x+iy$.
This inner product supplies an isomorphism between $S_2(N)$ and its dual:
for each linear map $L: S_2(N)\rightarrow\bo C$ there is a unique $\phi_L\in S_2(N)$  
such that $L(f)=\langle f,\phi_L\rangle $ for all $f\in S_2(N)$.  
Now let $F$ be a smooth function on~$\bo H$, not necessarily holomorphic, which transforms like a
modular form of weight~2 on $\Gamma_0(N)$ and which is `small' at the cusps of that group.
Below we will be more precise what we mean by `small'; it implies that the integral  \eqn1 with $g$ replaced by
 $F $ exists for all $f\in S_2(N)$, i.e. $\langle -,F\rangle$ is a linear functional on $S_2(N)$. 
The {\sl holomorphic projection} of $F$ is defined as the unique ${\tilde F}\in S_2(N)$ such that
$\langle f,F\rangle =\langle f,{\tilde F}\rangle$ for all $f\in S_2(N)$.
The Fourier coefficients of $\tilde F$ can be expressed in terms of the Fourier coefficients of $F$ and
the growth rate of $F$ at the cusps. 
It will turn out that we are only interested in the first Fourier coefficient of $F$, hence we restrict to that case. 
\bigskip\noindent\bf Theorem\,\thmnumber\sl 
Let  $N$ be a positive integer and let $F(z)=\sum_{m=-\infty}^\infty a_m(y)e^{2\pi imz}$ be a smooth function on
 $\bo H$ that transforms like a modular form of weight~$2$ on $\Gamma_0(N)$.
Suppose that for every positive divisor $M$ of $N$ there exist $\epsilon>0$ and complex numbers $A(M),B(M) $ such that 
$$
(cz+d)^{-2}F({az+b\over cz+d})=A(M)\log y+B(M)+O(y^{-\epsilon})\ \hbox{\ as }y=\Im(z)\rightarrow\infty\leqno(\eqnumber)
$$
for all $\ma abcd\in\SL_2(\bo Z)$ with $\gcd(c,N)=M$.
Let ${\tilde F}(z) =\sum_{m = 1} ^\infty a_me ^{2\pi imz} $ be the holomorphic projection of $F$. 
Then
$$
a_1\hskip -0.1 cm=\hskip -0.15cm\lim_{s\rightarrow 1\atop s\in\bo H_1}\hskip -0.15 cm
\Big[
\hskip-0.02cm 4\pi\hskip -0.1cm\int_0^\infty \hskip -0.35cm a_1(y)e^{-4\pi y}y^{s-1}dy
+{24\alpha\over s-1}
\Big]
-48\alpha\Big[
\sum_{p|N}{\log p\over p^2-1}+\log 2+{1\over 2}+
{\zeta^\prime\over\zeta}(2)\Big]+24\beta
$$
with $\bo H_1 =\{s\in\bo C:\Re(s) >1\} $ and
$$
\displaylines{
\textstyle\alpha=\prod\limits_{p|N}(1-p^{-2})^{-1}\sum\limits_{M|N}{\mu(M) A(M)\over M^2},\cr
\textstyle\beta=\prod\limits_{p|N}(1-p^{-2})^{-1}\sum\limits_{M|N}{\mu(M) [B(M)-2A(M)\log M]\over
M^2}.\cr
}$$
\noindent\bf Proof. \rm This is the holomorphic projection lemma for the first Fourier coefficient on page 534 of [\GZd]
(we corrected two typos).
Recall that a cusp form $f\in S_2(N)$ is exponentially small at the cusps.
Using \eqn1, which states that the growth rate of $F$ at the cusp  ${a\over c}\in\bo P ^1(\bo Q)$ only depends on
 $\gcd(c,N)$, we find that the $\Gamma_0(N)$-invariant function $f(z)\overline{F(z)}\Im(z)^2$
is bounded on $\bo H$. 
Therefore the inner product $\langle f,F\rangle$ converges and hence the holomorphic projection $\tilde F $ of $F$ is well defined. 

In order to sketch the proof we need the Poincar\'e series
$$
P(z)=\lim_{s\rightarrow1\atop s\in\bo H_1}
\sum_{\ma abcd\in\Gamma_\infty\back\Gamma_0(N)}
{1\over(cz+d)^2}{y^{s-1}\over|cz+d|^{2s-2}}e^{2\pi i{az+b\over cz+d}}
$$
where $\Gamma_\infty$ denotes the stabalizer of $\infty$ in $\Gamma_0(N)$.
This function is a holomorphic cusp form of weight 2 on $\Gamma_0(N)$. 
For each $f\in S_2(N)$ the inner product $\langle f,P\rangle$ is equal to $ 1\over 4\pi $ times the first Fourrier coefficient of $f$.
In particular we have $\langle{\tilde F},P\rangle={a_1\over4\pi}$.
If $F$ satisfies \eqn1 with $A(M)=B(M)=0$ for all $M$ then $\langle F,P\rangle$ is convergent
and equal to $\lim_{s\rightarrow1\atop s\in\bo H_1}\int_0^\infty a_1(y)e^{-4\pi y}y^{s-1}dy$.
In this special case both $\alpha$ and $\beta$ are zero and the theorem follows from the equality
$\langle {\tilde F},P\rangle=\langle F,P\rangle$.

In general we substract from $F$ a linear combination of non-holomorphic Eisenstein series of weight 2
to obtain a function $F^*$ which satisfies the asymptotic relation \eqn1 with $A(M)=B(M)=0$. 
As these Eisenstein series are perpendicular to the elements of $S_2(N)$, the functions $F$ and 
$F^*$ have the same holomorphic projection.     
The theorem follows by applying the arguments of the previous paragraph to the function $F^*$.
\hfill$\qed$

\medskip \noindent
To apply this theorem to $F(z) = F_{\go n}(z)$ we have to check condition \eqn1 for this function.
\bigskip\noindent\bf Proposition\,\thmnumber\sl 
Let $\go n$ be a primitive ideal of norm $N$ such that $\chi(\go n) =1$ and let $M$ be a positive divisor of $N$.
For all $\ma abcd\in\SL_2(\bo Z)$ with $(c,N)=M$ we have
$$\displaylines{
(cz+d)^{-2}F_\go n({az+b\over cz+d})=A(M)\log y+B(M)+O(y^{-1})\ \hbox{\ as }y=\Im(z)\rightarrow\infty\cr
\text{\sl with} \hfill \cr
\textstyle
A(M) =\varepsilon(M)M{A(N)\over N},\qquad
B(M) =\varepsilon (M)M\Big({B(N)\over N} +{A(N)\over N}\log{M\over N}\Big)\cr
\text{\sl  and}\hfill \cr
\textstyle
A(N) ={h_1 ^\prime h_2 ^\prime\over 2},\quad
B(N) ={h_1 ^\prime h_2 ^\prime\over 2}\Big(\half\log(DN) -\log\pi -\gamma +{L ^\prime\over L}(1,d_1) + 
{L ^\prime\over L}(1,d_2)\Big).\cr
} $$
\medskip \noindent\bf Proof. \rm 
The proposition follows from the Fourier expansions of $F_\go n(z)$  at the various cusps of the group $\Gamma_0(N)$.
Fix a matrix $A =\ma abcd\in\SL_2(\bo Z) $ and denote the greatest common divisor of $c$ and $N$ by $M$. 
It follows from the inclusion $\ma 1 {\go n}01 \subset A^{-1}\PSL_2(\ca O,\go n)A$
that the function 
$$
E_{\go n,s,A}(z_1,z_2)=(cz_1+d)^{-1}(cz_2+d)^{-1} E_{\go n,s}(Az_1,Az_2)
$$
satisfies $E_{\go n,s,A}(z_1+\lambda,z_2+\bar\lambda)=E_{\go n,s,A}(z_1,z_2)$
for all $\lambda\in\go n$.
It therefore has a Fourier expansion
$$
E_{\go n,s,A}(z_1,z_2)=\sum_{\nu\in (\go n\go d) ^{-1}} c_{\nu,s,A}(y_1,y_2)
e^{2\pi i (\nu x_1+\ba \nu x_2)}\leqno(\eqnumber)
$$
with 
$$
c_{\nu,s,A}(y_1,y_2) ={1\over N\sqrt{D}}\int_{\bo R ^2/\go n}E_{\go n,s,A}(z_1,z_2)e^{-2\pi i(\nu x_1+\ba \nu x_2)}dx_1dx_2\leqno(\eqnumber)
$$
for $\nu\in (\go n\go d)^{-1}$. 
The computation of this integral is very similar to the computation of  \eqn{16}  in the proof of proposition  \thmn4.
For $\Re(s) >1$  the function 
$
E_{\go n,s,A}(z_1,z_2)
$
has the same double sum representation as 
$E_{\go n,s}(z_1,z_2) $ but now the inner sum ranges over $(m, n)\in(\go a\times\go a)/\ca O ^* $ 
not both zero such that $dm -cn\in\go n\go a$.
In order to deal with this last condition we define $\go n_0$ as the greatest common divisor of the integral ideals 
$c\ca O$ and $\go n$ and let $\go n_1 =\go n\go n_0^{-1}$.
Because $\go n$ is primitive we have $N_{K/\bo Q}(\go n_0) =M$ and $\chi(\go n_0) =\varepsilon(M)$.  
As one easily proves for fixed $m_0\in \go n_0\go a/\go n\go a$ there is an unique 
$n_0\in \go a/\go n_1\go a $ such that
for $m\congr m_0\mod \go n\go a$ the congruence $dm\congr cn\mod \go n\go a $ is equivalent to
$n\congr n_0\mod \go n_1\go a $.
Using this equivalence and the techniques used in the proof of proposition  \thmn4 we find for  $c_{\nu,s,A}(y_1,y_2)$
$$
\eqalignno{
\textstyle\varepsilon(M) \big ( {M\over N} \big) ^{1+2s}L_K(1+2s,\chi)(y_1y_2) ^s +
{\textstyle\varepsilon(M)\Phi_s(0) ^2 L_K(2s,\chi)\over\textstyle NM^{2s-1}\sqrt{D}(y_1y_2)^{s} } &\ \ \ \text{\rm if\ }\nu =0;\cr
\noalign{\hbox{and}}
{\textstyle M\Phi_s(\nu y_1)\Phi_s(\ba \nu y_2)\over \textstyle N\sqrt{D}(y_1y_2)^{s} }
\sum_{[\go a]\in C} \sum_{m_0\in\go n_0/\go n}\hskip -0.4 cm\sum_{ {\scriptstyle m\in S\atop\scriptstyle  m\congr m_0(\go n)} \atop\scriptstyle \nu\in m(\go a\go n_1\go d)^{-1} }\hskip -0.4cm 
{\textstyle\chi(\go a)N(\go a)^{2s}e^{2\pi i\Tr ({\nu\over m} n_0)}\over\textstyle N(m) |N(m) |^{2s-1}} &\ \ \ \text{\rm if\ }\nu\ne 0,\cr 
}$$
where $S$ is a complete set of representatives for $(\go a -\{0\})/\ca O ^*$.
Note that for $\nu\ne 0$ the triple sum is a finite sum; the elements $m\in S $ such that $\nu\in m(\go a\go n_1\go d)^{-1}$
correspond bijectively to the integral ideals in ideal class of $\go a^{-1}$ that contain~$\nu\go n_1\go d$.

To find the Fourier expansion of 
$F_{\go n,A}(z) ={\sqrt{D}\over 8\pi ^2}{\partial\over\partial s}E_{\go n,s,A}(z,z) |_{s =0} $ we substitute $z_1 =z_2 $ in  \eqn2  and differentiate termwise.
It turns out that the growth rate of the function $F_{\go n,A}(z)$ for $y\rightarrow\infty$ is dominated by its zeroth
Fourier coefficient.
An elementary calculation using \eqn{12}  and \eqn9 yields
$$
\displaylines{
{\sqrt{D}\over 8\pi ^2}{\partial\over\partial s}c_{0,s,A}(y,y) |_{s =0} =A(M)\log y+B(M) \cr
\text{\rm  with}\hfill \cr
\textstyle A(M) ={\sqrt{D} M\varepsilon(M)\over 2\pi^2N}L_K(1,\chi)\cr
\text{\rm and}\hfill\cr
\textstyle B(M) ={\sqrt{D} M\varepsilon(M)\over 2\pi^2N}\Big( L_K ^\prime (1,\chi) + 
\big(\half\log(DN) -\log\big({N\over M}\big) -\log\pi -\gamma\big)L_K(1,\chi)\Big). \cr}
$$
By comparing the Euler product expansions we find the equality ${L(s,d_1)}{L(s,d_2)} =L_K(s,\chi) $, where the Dirichlet series $L(s,d_i)$ are defined as in proposition~\thmn8, and the  class number formula yields $L_K (1,\chi)={\pi ^2\over\sqrt{D}}h^\prime_1h ^\prime_2$.

To complete the proof it suffices to show
$$
\sumprime_{\nu\in(\go n\go d)^{-1}} {\partial\over\partial s}c_{\nu,s,A} |_{s =0} (y,y)e^{2\pi i x(\nu +\ba \nu )}=O(y^{-1})\hskip 0.6 cm\text{\rm  for\ }y\rightarrow\infty.\leqno(\eqnumber)
$$
By changing the path of integration and differentiating under the integral sign in~\eqn{15}  we obtain
$$
{\partial\over\partial s}\Phi_s(t) |_{s =0} =O(|t |^{-1}e^{-\pi |t |})\hskip1 cm\text{\rm  for all nonzero $t$}.
$$
(In fact one can replace  $\pi$ by any positive number smaller then $2\pi $.)
Using  \eqn{11} we find for all $y >1 $
$$
{\partial\over\partial s}c_{\nu,s,A} |_{s =0} (y,y)= O\Big(P(|\nu |, | {\ba \nu} |)e ^ { -y (|\nu | +|{\ba \nu} |)} \Big)
$$
for some polynomial $P\in\bo C[X,Y]$.
For $\nu={a +b\sqrt{D}\over2\sqrt{D}}\in (\go n\go d) ^{-1} $ we have $|\nu | + |{\ba \nu} | =\max({|a |\over\sqrt{D}}, |b |) $.
We conclude that $\#\{\nu\in (\go n\go d) ^{-1}:|\nu | + |{\ba \nu} | =x\} =O(x) $  for all $x\in\bo R_{>0} $ and \eqn1 follows easily.\hfill $\qed$
\bigskip\noindent\rm 
Before we state the main result of this section we introduce the following function.
For an ideal $\go a$ and $s\in\bo H_1=\{s\in\bo C:\Re(s) >1\}$ define 
$$
T_{\go a}(s) =\sum_{{\nu\in \go a\go d^{-1}\atop\nu >0 >\ba \nu}\atop\Tr(\nu) =1} \sigma_{0,\chi}\big ((\nu)\go d\go a ^{-1}\big)Q_{s-1}(1+2 |\ba \nu |)
$$
where $Q_{s-1}(t)$ is the Legendre function of the second kind defined by  \eqn{50} in section~3.
For $\nu$ as in the above sum, the pair $(\nu,\ba \nu) $ ranges over those elements in the lattice $\go a\go d^{-1}\subset\bo R ^2 $
that are on the half line $\{\Tr(\nu) =1:\nu >0 >\ba \nu\} \subset\bo R ^2 $.
Because  $\sigma_{0,\chi}(\go a) =O(N(\go a) ^\delta) $ for any $\delta>0$  and $Q_{s-1}(t) =O(t^{-s}) $ for $t\rightarrow\infty$, the above sum converges absolutely and locally uniformly on $\bo H_1$. 
\bigskip\noindent\bf Theorem\,\thmnumber\sl 
Let $\go n$ be a primitive ideal prime to $\go d$ of norm $N$  such that $\chi(\go n) =1$.
The first Fourier coefficient $a_1$ of the holomorphic projection of $F_{\go n}(z)$ is equal to:
$$
a_1=S_{\go n} -\lim_{s\rightarrow1\atop s\in\bo H_1}\Big [T_{\go n}(s) +T_{\overline {\go n}}(s)-{24\alpha\over s-1}\Big]
+12\alpha \Big (2C+2\sum_{p |N}{\log p\over \varepsilon(p)p +1} -\log N\Big)
$$
with $C$ as in theorem  $\thmn8$ and
$$
\alpha ={h_1 ^\prime h_2 ^\prime\over2N} \prod_{p |N}{p\over p +\varepsilon(p)}\hskip 0.5 cm\text{\sl and}\hskip 0.5 cm
S_{\go n}=\sum_{{\nu\in\go n\go d^{-1}\atop\nu \gg 0}\atop\Tr(\nu) =1} {\sigma}^\prime_\chi\big( (\nu)\go d\go n^{-1}\big).
$$
\medskip\noindent\bf Proof. \rm First we use proposition \thmn2 to compute $\alpha$ and $\beta$ defined in theorem $\thmn3 $ for the function
 $F(z) =F_{\go n}(z) $.
An easy computation yields the following formulas:
$$\displaylines{
\textstyle\alpha=\prod\limits _{p|N}(1-p^{-2})^{-1}\cdot{A(N)\over N}\sum\limits_{M|N}{\mu(M) \varepsilon(M)\over M} =
{h_1 ^\prime h_2 ^\prime\over2N}\prod\limits_{p |N}{1- \varepsilon(p)p^{-1}\over1-p^{-2}}\cr
\textstyle\beta =\alpha\Big [ \half\log({D\over N}) -\log\pi -\gamma +{L ^\prime\over L}(1,d_1) + 
{L ^\prime\over L}(1,d_2)+\sum\limits_{p |N} {\varepsilon(p)\log p \over p-\varepsilon(p)} \Big].\cr
} $$
As $\go n$ is prime to $\go d$ we have $\varepsilon(p) ^2 = 1 $ for all primes $p |N$ and the above expression
for $\alpha$ is equal to the one quoted in the theorem. 

By equations  \eqn{9} and \eqn{8} we find that the first Fourier coefficient of $F_{\go n}(z)$ is equal to
$$
\displaylines{
a_1(y)=S_{\go n}-\half\hskip -0.2 cm\sum_{{\nu\in \go n\go d^{-1}\atop\nu >0 >\ba \nu}\atop\Tr(\nu) =1}\hskip -0.1 cm\sigma_{0,\chi}\big( (\nu)\go d\go n^{-1}\big)\Phi(|{\bar\nu}|y)
-\half\hskip -0.2 cm\sum_{{\nu\in {\overline{\go n}}\go d^{-1}\atop\nu >0 >\ba \nu} \atop\Tr(\nu) =1}\hskip -0.1 cm\sigma_{0,\chi}\big( (\nu)\go d {\overline{\go n}}^{-1}\big)\Phi(|\ba\nu |y),\quad\cr
}$$
with $S_{\go n}$ as in the statement of the theorem and where $\Phi(t)$ is defined by  \eqn{7}.
The two sums on the right are infinite sums and converges because of \eqn{6}.
For $s\in\bo H_1$ we use the integral representation $\Gamma(s) =\int_0 ^\infty e^{-y}y ^ {s-1}dy $
and obtain
$$
\displaylines{
4\pi\hskip -0.1cm\int_0^\infty \hskip -0.1 cm a_1(y)e^{-4\pi y}y^{s-1}dy=
{\Gamma(s)\over(4\pi)^{s-1}} S_{\go n}-\half\hskip -0.2 cm\sum_{{\nu\in \go n\go d^{-1}\atop\nu >0 >\ba \nu}\atop\Tr(\nu) =1} \hskip -0.1 cm\sigma_{0,\chi}\big ((\nu)\go d\go n^{-1}\big)\Psi_{s-1}(|\ba \nu |)\hfill \cr
\hfill -\half\hskip -0.2 cm\sum_{{\nu\in {\overline{\go n}}\go d^{-1}\atop\nu >0 >\ba \nu}\atop\Tr(\nu) =1} \hskip -0.1 cm\sigma_{0,\chi}\big ((\nu)\go d{\overline{\go n}}^{-1}\big)\Psi_{s-1}(|\ba \nu |),\cr
}$$
with
$$
\Psi_{s-1} (\lambda) =4\pi\int_0 ^\infty\Phi(\lambda y)e^{-4\pi y}y^{s-1}dy\hskip 0.5 cm\text{\rm for\ } \lambda >0,\Re(s) >0.
$$
For large $\lambda$ and $s$ close to 1 this function differs little from the Legendre function of the second kind $Q_{s-1}(1+2\lambda)$.
More precisely, using Taylor expansions one can prove [\GZe, page 218]
$$\textstyle
\Psi_{s-1} (\lambda) -{2\Gamma(2s)\over (4\pi)^{s-1}\Gamma(s+1)}Q_{s-1}(1+2\lambda ) =
\cases 
O(\lambda^{-s-1}) &\text{\rm  for  $\Re(s) > 0$ and $\lambda\rightarrow\infty$ };\cr
0 &\text{\rm  for $s =1 $ and $\lambda>0$ },\cr
\endcases
$$
where the implied constant is independent of $s$.
This estimate allows us to replace $\Phi_{s-1}$ by $Q_{s-1}$: 
$$\displaylines{\textstyle
\lim\limits_{s\rightarrow 1\atop s\in\bo H_1}\hskip-0.05cm\Big[4\pi\int\limits_0^\infty a_1(y)e^{-4\pi y}y^{s-1}dy+{24\alpha\over s-1}\Big]=
S_{\go n} - \lim\limits_{s\rightarrow1\atop s\in\bo H_1}\hskip-0.05cm\Big [T_{\go n}(s) +T_{\overline{\go n}}(s)
- {(4\pi)^{s-1}\Gamma(s+1)\over\Gamma(2s)}{24\alpha\over s-1}\Big]\cr
} $$
If we substitute this limit into the formula for the first Fourier coefficients of the holomorphic projection of
 $F_{\go n}(z)$ given in theorem \thmn3  and use the Taylor expansion
${(4\pi)^{s-1}\Gamma(s+1)\over \Gamma(2s)}=1+(\log(4\pi) +\gamma -1  )(s-1)+O\big((s-1) ^2\big)$
we obtain 
$$
\displaylines{
a_1=S_{\go n} -\lim_{s\rightarrow1\atop s\in\bo H_1}\Big [T_{\go n}(s) +T_{\overline {\go n}}(s)-{24\alpha\over s-1}\Big]
-24\alpha\cdot\Big[2+2 {\zeta^\prime\over\zeta}(2)
-\half\log\big({D\over N}\big) - {L ^\prime\over L}(1,d_1)\cr
\hfill  -{L ^\prime\over L}(1,d_2) +2\sum_{p|N}{\log p\over p^2-1} -\sum_{p |N} {\varepsilon(p)\log p \over p-\varepsilon(p)} \Big]\cr
} $$
which is equivalent to the statement in the theorem.\hfill $\qed$
  
\head 5. Conclusion of the proofs
\endhead 
\noindent 
Let $d_1$ and $d_2$ the two negative fundamental discriminants that are relatively prime and let $K$ be the real quadratic field of discriminant $D =d_1d_2 $.
By making the appropriate choice for the  $K$-ideal $\go n$  we can use theorem \thmn1 to calculate the limits
occurring on the right hand side of the equations in the theorems~\thmn8 and~\thmn6, which will complete the proofs of theorems~4 and~5 in the introduction.

Let $p$ be a rational prime and let  $a$ be a positive integer that satisfy the following properties:
\item{1.}  $p$ splits completely in $K/\bo Q$; Denote the primes above $p$ by $\go p$ and $\go q$;
\item{2.} $a$ is even in case $\varepsilon(p) = -1$;
\item{3.} the vector space $S_2(p ^a)$ of holomorphic modular forms of weight 2 on $\Gamma_0(p ^a)$ 
is zero dimensional. 

\noindent Recall that property 3 is equivalent with the condition that the compact Riemann surface
 $\Gamma_0(p ^a)\back\bH $ has genus zero and that there are only finitely many pairs $(p,a)$ with this property.
Because of the third condition on $p^a$, the first Fourier coefficient of the holomorphic projection of $F_{\go p ^a}$ and $F_{\go q ^a}$  are both zero.
If we add the result of theorem \thmn1 for the ideals $\go p ^a $ and $\go q ^a $ we obtain
$$
2\lim_{s\rightarrow1\atop s\in\bo H_1}\Big [T_{\go p ^a}(s) +T_{\go q ^a}(s)-{24\alpha\over s-1}\Big] =
S_{\go p ^a} +S_{\go q ^a}
+ 24\alpha \Big (2C+2{\log p\over \varepsilon(p)p +1} -\log p ^a\Big)
$$
with $C$ as in theorem \thmn8 and
$\alpha ={h_1^\prime h_2^\prime \over 2p ^ {a -1} (p +\varepsilon(p))}$,
$$
T_{\go p ^a}\hskip-0.05cm(s) +T_{\go q ^a }\hskip-0.05cm(s)
\hskip-0.05cm =
\hskip -0.4 cm\sum_{{\nu\in \go p ^ a\go d^{-1}\atop\nu >0 >\ba \nu}\atop\Tr(\nu) =1}
\hskip -0.35 cm \sigma_{0,\chi}\big ((\nu)\go d\go p^{-a}\big)Q_{s-1}\hskip-0.05cm(1+2 |\ba \nu |)
\hskip0.05cm+\hskip -0.35 cm\sum_{{\nu\in \go q ^ a\go d^{-1}\atop\nu >0 >\ba \nu}\atop\Tr(\nu) =1}
\hskip -0.35 cm \sigma_{0,\chi}\big ((\nu)\go d\go q^{-a}\big)Q_{s-1}\hskip-0.05cm(1+2 |\ba \nu |)
$$
and
$$
S_{\go p ^a} +S_{\go q ^a} =
\hskip -0.1 cm\sum_{{\nu\in\go p ^a\go d^{-1}\atop\nu \gg 0}\atop\Tr(\nu) =1} 
\hskip -0.2 cm {\sigma}^\prime_\chi\big( (\nu)\go d\go p^{-a}\big)
+\hskip -0.1 cm\sum_{{\nu\in\go q ^a\go d^{-1}\atop\nu \gg 0}\atop\Tr(\nu) =1}
\hskip -0.2 cm {\sigma}^\prime_\chi\big( (\nu)\go d\go q^{-a}\big).
$$
The totally positive elements $\nu\in\go d^{-1}$ of trace 1 are of the form $\nu={n+\sqrt{D}\over2\sqrt{D}}$ 
with~$n$ a rational integer satisfying $| n | <\sqrt{D} $ and $n\congr D\mod 2$ .
Such an element $\nu $ is not divisible by an integer, hence its norm is divisible by $p ^a $ if and only if
either $\nu\in\go p ^a\go d^{-1}$ or $\nu\in\go q ^a\go d ^{-1} $.
For $\nu={n+\sqrt{D}\over2\sqrt{D}}\in\go p ^a\go d^{-1} $  the norm maps the ideals dividing the primitive integral ideal $(\nu)\go d\go p ^{-a} $ bijectively to the positive divisors of $N_{K/\bo Q}\big((\nu)\go d\go p^{-a}\big) ={D-n ^2\over 4p ^a} $.
Hence we find for these $\nu$ 
$$
{\sigma} ^\prime_\chi\big ((\nu)\go d\go p^{-a} \big) =\sum_{\go a |(\nu)\go d\go p^{-a}}\chi(\go a)\log N_{K/\bo Q}(\go a) =\sum_{d |{D-n ^2\over 4p ^a}}\varepsilon(d)\log d.
$$ 
The same formula holds for $\sigma ^\prime_\chi\big((\nu)\go d\go q^{-a}\big)$ in case $\nu\in\go q^a\go d^{-1} $ 
and we obtain
$$
S_{\go p ^a} +S_{\go q ^a} =\sum_{|n| <\sqrt{D}}\sum_{d |{D-n ^2\over 4p ^a}}\varepsilon(d)\log d.
$$
In a similar way we find
$$
T_{\go p ^a}(s) +T_{\go q ^a}(s) =\sum_{n >\sqrt{D}}\sum_{d |{n ^2-D\over 4p ^a}}\varepsilon(d)
Q_{s-1}\Big({n\over\sqrt{D}}\Big)  
$$
and the above limit can be rewritten as
$$\eqalign{
\qquad\lim\limits_{s\rightarrow1\atop s\in\bo H_1}
\hskip-0.0cm\left [
{48\alpha\over s-1} -2
\hskip-0.1cm\sum\limits_{n >\sqrt{D}}
\hskip -0.0 cm\rho_{p ^a}(n) Q_{s-1}\Big({n\over\sqrt{D}}\Big) \right] 
=\hfill&\cr
&\hskip-4.9cm 24\alpha \Big (\log p ^a- {2\log p\over \varepsilon(p)p +1} -2C\Big)-
\hskip -0.2 cm\sum\limits_{|n|<\sqrt{D}}
\sum\limits_{d |{D-n ^2\over 4p ^a}}\varepsilon(d)\log d
}\leqno(\eqnumber)
$$
with 
$$
\rho_{p ^a}(n) =
\cases
\sum_{d |{n ^2-D\over4p ^a}}\varepsilon(d) &\text{\rm  if\ } n ^2\congr D\mod 4p ^a\cr
0 &\text{\rm  otherwise.} 
\endcases
$$

\noindent Now assume that both $d_1$ and $d_2$ are congruent to 2 modulo 3.
In this case the rational prime $3$ splits completely in $K/\bo Q$ and $\varepsilon(3) = -1 $.
As $\Gamma_0(9)\back\bH $ has genus~0 [\MI, \S 4.2] the pair 
$(p,a) =(3, 2) $ satisfies the three properties that we stated in the begin of this section.
If we apply \eqn1 for this pair and use the equalities $\alpha={h_1^\prime h_2^\prime\over 12}$
and $\rho_{9}(n)=\rho^{\gamma_2}(n)$ (see \eqn{29}) we arrive at
$$\textstyle
\lim\limits_{s\rightarrow1\atop s\in\bo H_1}
\hskip-0.15cm\Big [
{4h_1^\prime h_2^\prime \over s-1} -2
\hskip-0.2cm\sum\limits_{n >\sqrt{D}}
\hskip-0.2cm\rho^{\gamma_2}(n) Q_{s-1}\hskip-0.05cm\big ({n\over\sqrt{D}}\big) 
\Big] 
\hskip-0.03cm=\hskip-0.03cm
2h_1^\prime h_2^\prime\big(3\log 3-2C\big)-
\hskip-0.2cm\sum\limits_{|n| <\sqrt{D}}\sum\limits_{d |{D-n ^2\over 36}}
\hskip-0.25cm\varepsilon(d)\log d.
$$
Substituting this limit into the formula of theorem \thmn8 yields
$$\eqalign{
\log \N\big(\gamma_2(\alpha_1),\gamma_2(\alpha_2)\big)^{8\over w_1w_2}
&=6h_1^\prime h_2^\prime\log 3-
\hskip-0.2cm\sum\limits_{|n| <\sqrt{D}}\sum\limits_{d |{D-n ^2\over 36}}
\hskip-0.2cm\varepsilon(d)\log d\cr
&=6h_1^\prime h_2^\prime\log 3+
\hskip-0.2cm\sum\limits_{ n ^2 <{D}}\log F\Big( {D-n ^2\over 36} \Big)\cr
}
$$
where the last line follows from the equality $\varepsilon({D-n ^2\over36}) = -1 $ [\CO, page 306].
This concludes the proof of theorem 4.

Next we assume that both $d_1$ and $d_2$ are congruent to 1 modulo 8,
in particular $h_i ^\prime =h_i$.
Both pairs $(p,a) =(2,1) $ and $(p,a) =(2,2) $ satisfy the desired properties.
For $(p,a) =(2,1) $ we find $\alpha={h_1 h_2\over 6}$ and because of the equality $2\rho_2(n)=\rho^{\omega}(n)$ (see \eqn{26}) 
equation \eqn1 now reads
$$\textstyle
\lim\limits_{s\rightarrow1\atop s\in\bo H_1}
\hskip-0.15cm\Big [
{8h_1h_2\over s-1} -
\hskip-0.2cm\sum\limits_{n >\sqrt{D}}
\hskip-0.2cm\rho^{\omega}(n) Q_{s-1}\big ({n\over\sqrt{D}}\big) \Big] 
=
4h_1h_2\big({1\over 3}\log 2-2C\big)-
\hskip-0.03cm \hskip-0.3cm\sum\limits_{|n |<\sqrt{D}}\sum\limits_{d |{D-n ^2\over 8}}
\hskip-0.2cm\varepsilon(d)\log d.
$$
Applying \eqn1 for $(p,a) =(2, 2)$, $\alpha={h_1 h_2\over 12}$ and
$\rho_4(n)=\rho^{\omega_2}(n)$ (see \eqn{24}) and dividing by 2 yields
$$\textstyle
\lim\limits_{s\rightarrow1\atop s\in\bo H_1}
\hskip-0.15cm\Big [
{2h_1 h_2\over s-1} -
\hskip-0.2cm\sum\limits_{n >\sqrt{D}}
\hskip-0.2cm\rho^{\omega_2}(n) Q_{s-1}\big ({n\over\sqrt{D}}\big) \Big] 
=
2h_1h_2\big({2\over 3}\log 2-C\big) -
\hskip-0.2cm\sum\limits_{0<n<\sqrt{D}}\sum\limits_{d |{D-n ^2\over 16}}
\hskip-0.2cm\varepsilon(d)\log d.
$$
If we substitute these limits in the formulas of theorem  \thmn6 we find
$$\eqalignno{
\log \N\big(\omega(\alpha_2), \omega(\alpha_2)\big)&=
12h_1h_2\log 2-
\hskip-0.2cm \sum\limits_{|n |<\sqrt{D}}\sum\limits_{d |{D-n ^2\over 8}}
\hskip-0.2cm\varepsilon(d)\log d\cr
\noalign{\hbox{\text{\rm and}}}
\log \N\big(\omega_2(\alpha_2), \omega_2(\alpha_2)\big)&=
-
\hskip-0.2cm\sum\limits_{0 <n <\sqrt{D}}\sum\limits_{d |{D-n ^2\over 16}}
\hskip-0.2cm\varepsilon(d)\log d,\cr
}$$
from which theorem 5 follows.

\vfill\eject

\enddocument